\documentstyle[12pt]{article}
\renewcommand{\baselinestretch}{1.3}

\newtheorem {th}{Theorem}[section]
\newtheorem {lem}[th]{Lemma}

\newtheorem{conj}{Conjecture}
\newcounter{pfigure}
\def\lap{\bigtriangleup}

\def\VV{{\overline V}}
\def\WW{{\overline W}}
\def\xy{{\vec {xy}}}
\def\xz{{\vec {xz}}}
\def\yx{{\vec {yx}}}
\def\zw{{\vec {zw}}}

\def\ww{\vec {w}}
\def\Cox{\hfill \Box}

\def\dconv{\, {\stackrel {{\cal D}} {\rightarrow}}}

\def\dd{\delta}
\def\ee{\epsilon}
\def\E{{\bf{E}}}
\def\P{{\bf{P}}}
\def\pois{{\cal P}_1}
\def\N{\hbox{I\kern-.2em\hbox{N}}}
\def\R{\hbox{I\kern-.2em\hbox{R}}}
\def\Z{\hbox{Z\kern-.4em\hbox{Z}}}
\def\A{{\cal{A}}}
\def\C{{\cal{C}}}

\def\|{\, | \, }

\def\Tree{{\bf T}}

\begin{document}

\section{Introduction} \label{1}

There are several good reasons you might want to read about uniform
spanning trees, one being that spanning trees are 
useful combinatorial objects.  Not only are they fundamental in
algebraic graph theory and combinatorial geometry, but they 
predate both of these subjects, having been used by Kirchoff in
the study of resistor networks.  This article addresses the question
about spanning trees most natural to anyone in probability 
theory, namely what does a typical spanning tree look like?

Some readers will be happy to know that understanding the basic
questions requires no background knowledge or technical expertise.
While the model is elementary, the answers are surprisingly rich.
The combination of a simple question and a structurally complex
answer is sometimes taken to be the quintessential mathematical
experience.  This nonwithstanding, I think the best reason to
set out on a mathematical odyssey is to enjoy the ride.
Answering the basic questions about spanning trees depends on a sort of 
vertical integration of techniques and results from diverse realms of 
probability theory and discrete mathematics.  Some of the topics
encountered en route are random walks, resistor networks,
discrete harmonic analysis, stationary Markov chains, circulant
matrices, inclusion-exclusion, branching processes and the method of
moments.  Also touched on are characters of abelian groups, entropy
and the infamous incipient infinite cluster.  

The introductory section defines the model and previews some of the
connections to these other topics.  The remaining sections develop
these at length.  Explanations of jargon and results borrowed from other
fields are provided whenever possible.  Complete proofs are given in
most cases, as appropriate.

\subsection{Defining the model} \label{1.1}

Begin with a finite graph $G$.  That means a finite collection $V(G)$ 
of {\em vertices} along with a finite collection $E(G)$ of {\em edges}.
Each edge either connects two vertices $v$ and $w \in V(G)$ or else
is a {\em self-edge}, connecting some $v \in V(G)$ to itself.  There
may be more than one edge connecting a pair of vertices.  Edges
are said to be {\em incident} to the vertices they connect.  To make 
the notation less cumbersome we will write $v \in G$ and $e \in G$
instead of $v \in V(G)$ and $e \in E(G)$.  For $v,w \in G$ say $v$
is a {\em neighbor} of $w$, written $v \sim w$ if and only if some edge connects
$v$ and $w$.  Here is an example of a graph $G_1$ which will serve often as
an illustration.

\setlength{\unitlength}{2pt}
\begin{picture}(170,60)
\put(100,50){\circle{2}}
\put(98,53){A}
\put(140,50){\circle{2}}
\put(140,53){B}
\put(140,10){\circle{2}}
\put(140,3){C}
\put(100,10){\circle{2}}
\put(98,2){D}
\put(70,30){\circle{2}}
\put(62,28){E}
\put(100,50){\line(1,0){40}}
\put(120,53){$e_1$}
\put(100,50){\line(0,-1){40}}
\put(102,30){$e_6$}
\put(140,50){\line(0,-1){40}}
\put(142,30){$e_2$}
\put(140,10){\line(-1,0){40}}
\put(120,1){$e_3$}
\put(100,50){\line(-3,-2){30}}
\put(77,42){$e_5$}
\put(100,10){\line(-3,2){30}}
\put(77,12){$e_4$}
\end{picture}
figure~\thepfigure
\label{pfig1.1}
\addtocounter{pfigure}{1}

\noindent{Its} vertex set is $\{ A,B,C,D,E \}$ and it has six edges
$e_1 , \ldots , e_6$, none of which is a self-edge.  

A subgraph of a graph $G$ will mean a graph with the same vertex
set but only a subset of the edges.  (This differs from standard
usage which allows the vertex set to be a subset as well.)
Since $G_1$ has 6 edges, there are $2^6 = 64$ possible different
subgraphs of $G_1$.  A subgraph $H \subseteq G$
is said to be a {\em forest} if there are no cycles, i.e. you cannot
find a sequence of vertices $v_1 , \ldots , v_k$ for which there
are edges in $H$ connecting $v_i$ to $v_{i+1}$ for each $i < k$ and an 
edge connecting $v_k$ to $v_1$.  In particular $(k=1)$ there are no
self-edges in a forest.  A {\em tree} is a forest that is connected, 
i.e. for any $v$ and $w$ there is a path of edges that connects them.
The {\em components} of a graph are the maximal connected subgraphs,
so for example the components of a forest are trees.
A {\em spanning} forest is a forest in which every vertex has at
least one incident edge; a {\em spanning tree} is a tree in which
every vertex has at least one incident edge.  If $G$ is connected
(and all our graphs will be) then a spanning tree is just a subgraph
with no cycles such that the addition of any other edge would create
a cycle.  From this it is easy to see that every connected graph 
has at least one spanning tree.   

Now if $G$ is any finite connected graph, imagine listing all of its
spanning trees (there are only finitely many) and then choosing one
of them at random with an equal probability of choosing any one.
Call this random choice $\Tree$ and say that
$\Tree$ is a {\em uniform random spanning tree} for $G$.
In the above example there are eleven spanning trees for $G_1$
given (in the obvious notation) as follows:
$$ \begin{array} {rrrr}
e_1e_2e_3e_4 & e_1e_2e_3e_5 & e_1e_2e_4e_5 & e_1e_3e_4e_5 \\[2ex]
e_2e_3e_4e_5 & e_1e_2e_4e_6 & e_1e_3e_4e_6 & e_2e_3e_4e_6 \\[2ex]
e_1e_2e_5e_6 & e_1e_3e_5e_6 & e_2e_3e_5e_6 
\end{array} $$
In this case, $\Tree$ is just one of these eleven trees, picked
with uniform probability.  The model is so simple, you may wonder what
there is to say about it!  One answer is that the model has some 
properties that are easy to state but hard to prove; these are
introduced in the coming subsections.  Another answer is that the
definition of a uniform random spanning tree does not give us
a way of readily computing {local characteristics} of the 
random tree.  To phrase this as a question: can you compute probabilities
of events local to a small set of edges, such as $\P (e_1 \in \Tree)$
or $\P (e_1 , e_4 \in \Tree)$ without actually enumerating all
of the spanning trees of $G$?  In a sense, most of the article is
devoted to answering this question.  (Events such as $e_1$ being in the
tree are called local in contrast to a global event such as the 
tree having diameter -- longest path between two vertices -- at most three.)

\subsection{Uniform spanning trees have negative correlations} \label{1.2}

Continuing the example in figure~1,
suppose I calculate the probability that $e_1 \in \Tree$.  That's easy:
there are 8 spanning trees containing $e_1$, so
$$\P (e_1 \in \Tree) = {8 \over 11}.$$
Similarly there are 7 spanning trees containing $e_4$ so
$$\P (e_4 \in \Tree) = {7 \over 11}.$$
There are only 4 spanning trees containing both $e_1$ and $e_4$, so 
$$ \P (e_1 \in \Tree \mbox{ and } e_4 \in \Tree) = {4 \over 11} .$$ 
Compare the probability of both of these edges being in the tree with
the product of the probabilities of each of the edges being in the tree:
$${8 \over 11} \cdot {7 \over 11} = {56 / 121} > {4 \over 11} . $$
Thus
$$\P (e_1 \in \Tree \| e_4 \in \Tree) = { \P (e_1 \in \Tree \mbox{ and }
e_4 \in \Tree) \over \P (e_4 \in \Tree)} < \P (e_1 \in \Tree ) $$
or in words, the conditional probability of $e_1$ being in the tree if you 
know that $e_4$ is in the tree is less than the original unconditional
probability.  This negative correlation of edges holds in 
general, with the inequality not necessarily strict.

\begin{th} \label{neg cor}
For any finite connected graph $G$, let $\Tree$ be a uniform spanning
tree.  If $e$ and $f$ are distinct edges, then $\P (e,f \in \Tree ) 
\leq \P (e \in \Tree) \P (f \in \Tree)$.
\end{th}

Any spanning tree of an $n$-vertx graph contains $n-1$ edges, so
it should seem intuitively plausible -- even obvious -- that if
one edge is forced to be in the tree then any other edge is less
likely to be needed.  Two proofs will be given later, but neither
is straightforward, and in fact the only proofs I know involve
elaborate connections between spanning trees, random walks and
electrical networks.  Sections~\ref{2} and~\ref{3} will be occupied with the
elucidation of these connections.  The connection between random walks 
and electrical networks will be given more briefly, since an excellent 
treatment is available \cite{DS}.  

As an indication that the previous theorem is not trivial, here is
a slightly stronger statement, the truth or falsity of which is
unknown.  Think of the distribution of $\Tree$ as 
a probability distribution on the outcome space $\Omega$ consisting
of all the $2^{|E(G)|}$ subgraphs
of $G$ that just happens to give probability zero to any 
subgraph that is not a spanning tree.  An event $A$ (i.e. any subset of
the outcome space) is called an {\em up-event} -- short for {\em
upwardly closed} -- if whenever a subgraph $H$ of $G$ has
a further subgraph $K$ and $K \in A$, then $H \in A$.
An example of an up-event is the event of containing at least two of the
three edges $e_1, e_3$ and $e_5$.  Say an event $A$ ignores an edge
$e$ if for every $H$, $H \in A \Leftrightarrow H \cup 
e \in A$.

\begin{conj} \label{stoch decr}
For any finite connected graph $G$, let $\Tree$ be a uniform spanning
tree.  Let $e$ be any edge and $A$ be any up-event that ignores $e$.
Then 
$$ \P (A \mbox{ and } e \in \Tree) \leq \P (A) \P (e \in \Tree) . $$
\end{conj}
Theorem~\ref{neg cor} is a special case of this when $A$ is the
event of $f$ being in the tree.  The conjecture is known to be true for 
{\em series-parallel} graphs and it is also know to be true in 
the case when $A$ is an {\em elementary cylinder event}, i.e. the
event of containing some fixed $e_1 , \ldots , e_k$.  On the negative
side, there are natural generalizations of graphs and spanning trees, namely
{\em matroids} and {\em bases} (see \cite{Wh} for definitions), 
and both Theorem~\ref{neg cor} and Conjecture~1 fail to generalize
to this setting.  If you're interested in seeing the counterexample,
look at the end of \cite{SW}.

\subsection{The transfer-impedance matrix} \label{1.3}

The next two paragraphs discuss a theorem that computes probabilities
such as $\P (e,f \in \Tree)$.  These computations alone would render
the theorem useful, but it appears even more powerful in the context
of how strongly it constrains the probability measure governing
$\Tree$.  Let me elaborate.

Fix a subset $S = \{ e_1 , \ldots , e_k \}$ of the edges of a 
finite connected graph $G$.  If $\Tree$ is a uniform random spanning
tree of $G$ then the knowledge of whether $e_i \in \Tree$ for each $i$
partitions the space into $2^k$ possible outcomes.  (Some of these
may have probability zero if $S$ contain cycles,
but if not, all $2^k$ may be possible.)  In any case, choosing
$\Tree$ from the uniform distribution on spanning trees of $G$
induces a probability distribution on $\Omega$, the space of
these $2^k$ outcomes.  There
are many possible probability distributions on $\Omega$:
the ways of choosing $2^k$ nonnegative numbers summing
to one are a $2^k - 1$-dimensional space.  Theorem~\ref{neg cor}
shows that the actual measure induced by $\Tree$ satisfies 
certain inequalities, so not all probability distributions on $\Omega$
can be gotten in this way.  But the set of probability distributions
on $\Omega$ satisfying these inequalities is still $2^k -1$-dimensional.
It turns out, however, that the set of probability distributions on $\Omega$
that arise as induced distributions of uniform spanning trees on
subsets of $k$ edges actually has at most the much smaller dimension 
$k(k+1)/2$.  This is a consequence of the following theorem which 
is the bulwark of our entire discussion of spanning trees:

\begin{th}[Transfer-Impedance Theorem] \label{exist matrix}
Let $G$ be any finite connected graph.  There is a symmetric function
$H(e,f)$ on pairs of edges in $G$ such that for any $e_1 , \ldots , e_r
\in G$, 
$$ \P (e_1 , \ldots , e_r \in \Tree) = \det M(e_1 , \ldots , e_r)$$
where $M(e_1 , \ldots , e_r)$ is the $r$ by $r$ matrix whose
$i,j$-entry is $H(e_i , e_j)$.  
\end{th}
By inclusion-exclusion, the probability of any event in $\Omega$ 
may be determined from the probabilities of $\P (e_{j_1} , \ldots , 
e_{j_r} \in \Tree)$ as $e_{j_1} , \ldots , e_{j_r}$ vary over all
subsets of $e_1 , \ldots , e_k$.  The theorem says that these
are all determined by the $k(k+1)/2$ numbers $\{ H(e_i , e_j)  :
i,j \leq k \}$, which shows that there are indeed only $k(k+1)/2$
degrees of freedom in determining the measure on $\Omega$.  

Another way of saying this is that the measure is almost completely 
determined by its
two-dimensional marginals, i.e. from the values of $\P (e,f \in \Tree)$
as $e$ and $f$ vary over pairs of (not necessarily distinct) edges.
To see this, calculate the values of $H(e,f)$.  The values of
$H(e,e)$ in the theorem must be equal 
to $\P (e \in \Tree)$ since $\P (e,e) = \det M(e) = H(e,e)$.  
To see what $H(e,f)$ is for $e \neq f$, write
\begin{eqnarray*}
\P (e,f \in \Tree) & = & \det M(e,f) \\[2ex]
 & = & H(e,e) H(e,f) - H(e,f)^2 \\[2ex]
& = & \P(e \in \Tree) \P (f \in \Tree) - H(e,f)^2
\end{eqnarray*}
and hence
$$H(e,f) = \pm \sqrt{\P (e \in \Tree) \P (f \in \Tree) - \P (e,f \in \Tree)}.$$
Thus the two dimension marginals determine $H$ up to sign, and $H$
determines the measure.  Note that the above square
root is always real, since by Theorem~\ref{neg cor} the quantity
under the radical is nonnegative.  
Section~\ref{4} will be devoted to proving Theorem\ref{exist matrix}, 
the proof depending heavily on the connections to random walks and 
electrical networks developed in Sections~\ref{2} and~\ref{3}.

\subsection{Applications of transfer-impedance to limit theorems} \label{1.4}

Let $K_n$ denote the complete graph on $n$ vertices, i.e. there
are no self-edges and precisely one edge connecting each pair of 
distinct vertices.  Imagine picking a uniform random spanning
tree of $K_n$ and letting $n$ grow to infinity.  What kind of
limit theorem might we expect?  Since a spanning tree of $K_n$ has only
$n-1$ edges, each of the $n(n-1)/2$ edges should have probability
$2/n$ of being in the tree (by symmetry) and is hence decreasingly
likely to be included as $n \rightarrow \infty$.  
On the other hand, the number of edges
incident to each vertex is increasing.  Say we fix a particular
vertex $v_n$ in each $K_n$ and look at the number of edges incident
to $v_n$ that are included in the tree.  Each of $n-1$ incident
edges has probability $2/n$ of being included, so the expected number of
of such edges is $2(n-1)/n$, which is evidently converging to 2.  If the
inclusion of each of these $n-1$ edges in the tree were independent
of each other, then the number of edges incident to $v_n$ in $\Tree$ 
would be a binomial random variable with parameters $(n-1 , 2/n)$;
the well known Poisson limit theorem would then say
that the random variable $D_\Tree (v_n)$ counting how many edges incident
to $v_n$ are in $\Tree$ converged as $n \rightarrow \infty$ to 
a Poisson distribution with mean two.  
(A quick explanation: integer-valued random variables $X_n$ are said to
converge to $X$ in distribution if $\P (X_n = k) \rightarrow \P (X=k)$
for all integers $k$.  In this instance, convergence of $D_\Tree (v_n)$
to a Poisson of mean two would mean that for each $k$,
$\P (D_\Tree (v_n) = k) \rightarrow e^{-2} k^2 / 2$ as $n \rightarrow
\infty$ for each integer $k$.)  Unfortunately
this can't be true because a Poisson(2) is sometimes zero, whereas
$D_\Tree (v_n)$ can never be zero.  It has however been shown \cite{Al}
that $D_\Tree (v_n)$ converges in distribution to the next simplest
thing: one plus a Poisson of mean one.  

To show you why this really
is the next best thing, let me point out a property of the
mean one Poisson distribution.  Pretend that if you picked a family
in the United States at random, then the number of children in the
family would have a Poisson distribution with mean one (population 
control having apparently succeeded).  Now
imagine picking a child at random instead of picking a family
at random, and asking how many children in the family.  You would
certainly get a different distribution, since you couldn't ever
get the answer zero.  In fact you would get one plus a Poisson
of mean one.  (Poisson distributions are the only ones with
this property.)  Thus a Poisson-plus-one distribution is a more
natural distribution than it looks at first.
At any rate, the convergence theorem is 

\begin{th} \label{Al Pois}
Let $D_\Tree (v_n)$ be the random degree of the vertex $v_n$ in a 
uniform spanning tree of $K_n$.  Then as $n \rightarrow \infty$,
$D_\Tree (v_n)$ converges in distribution to $X$ where $X$ is one
plus a Poisson of mean one.
\end{th}

Consider now the $n$-cube $B_n$.  Its vertices are defined to be
all strings of zeros and ones of length $n$, where two vertices are
connected by an edge if and only if they differ in precisely one
location.  Fix a vertex $v_n \in B_n$ and play the same game: choose a
uniform random spanning tree and let $D_\Tree (v_n)$ be the random 
degree of $v_n$ in the tree.  It is not hard to see again that the
expected value, $\E D$, converges to 2 as $n \rightarrow \infty$.
Indeed, for any graph the number of vertices in a spanning tree
is one less than the number of vertices, and since each edge
has two endpoints the average degree of the vertices will be
$\approx 2$; if the graph is symmetric, each vertex will then have the
same expected degree which must be 2.  One could expect 
Theorem~\ref{Al Pois} to hold for $B_n$ as well as $K_n$ and in fact it
does.  A proof of this for a class of sequences of graphs that includes
both $K_n$ and $B_n$ and does not use transfer-impedances 
appears in \cite{Al} along with the conjecture that the result
should hold for more general sequences of graphs.  This can indeed
be established, and in Section~\ref{5} we will discuss the proof of
Theorem~\ref{Al Pois} via transfer-impedances which can
be extended to more general sequences of graphs. 

The convergence in distribution of $D_\Tree (v_n)$ in these
theorems is actually a special case of a stronger kind of
convergence.  To begin discussing this stronger kind of convergence, 
imagine that we pick a uniform random spanning tree of a graph, say
$K_n$, and want to write down what it looks like ``near $v_n$''.
Interpret ``near $v_n$'' to mean within a distance of $r$ of 
$v_n$, where $r$ is some arbitrary positive integer.
The answer will be a rooted tree of height $r$.  (A
{\em rooted} tree is a tree plus a choice of one of its vertices, called
the root.  The {\em height} of a rooted tree is the maximum distance
of any vertex from the root.)  The rooted tree representing $\Tree$
near $v_n$ will be the tree you get by picking up $\Tree$, dangling
it from $v_n$, and ignoring everything more than $r$ levels
below the top. 

Call this the {\em $r$-truncation} of $\Tree$, written $\Tree 
\wedge_{v_n} r$ or just $\Tree \wedge r$ when the choice of $v_n$ 
is obvious.  For example, suppose $r=2$, $v_n$ has 2 neighbors in 
$\Tree$, $w_1$ and $w_2$, $w_1$ has 3 neighbors other than $v_n$ 
in $\Tree$ and $w_2$ has none.  This information is encoded
in the following picture.  The picture could also have been
drawn with left and right reversed, since we consider this to
be the same abstract tree, no matter how it is drawn. 

\setlength{\unitlength}{1pt}
\begin{picture}(310,120)(-90,-10)
\put(90,65) {\circle*{5}}
\put(95,68) {$v_n$}
\put(90,65) {\line(2,-1){60}}
\put(90,65) {\line(-2,-1){60}}
\put(150,35) {\line(1,-1){30}}
\put(150,35) {\line(0,-1){30}}
\put(150,35) {\line(-1,-1){30}}
\put(150,35) {\circle{2}}
\put(30,35) {\circle{2}}
\put(180,5) {\circle{2}}
\put(150,5) {\circle{2}}
\put(120,5) {\circle{2}}
\end{picture}
\setlength{\unitlength}{2pt}
figure~\thepfigure
\label{pfig1.2}
\addtocounter{pfigure}{1}

When $r=1$, the only information in $\Tree \wedge r$ is the number
of children of the root, i.e. $D_\Tree (v_n)$.  Thus the
previous theorem asserts the convergence in distribution of
$\Tree \wedge_{v_n} 1$ to a root with a (1+Poisson) number of
vertices.  Generalizing this is the following theorem, proved in
Section~\ref{5}.

\begin{th} \label{pois conv}
For any $r \geq 1$, as $n \rightarrow \infty$, $\Tree \wedge_{v_n} r$
converges in distribution to a particular random tree, $\pois
\wedge r$ to be defined later.
\end{th}
Convergence in distribution means that for any fixed tree $t$ of height
at most $r$, $\P (\Tree \wedge_{v_n} r = t)$ converges as $n \rightarrow
\infty$ to the probability of the random tree $\pois \wedge r$ equalling
$t$.  As the notation indicates, the random tree $\pois \wedge r$ is
the $r$-truncation of an infinite random tree.  It is in fact the
tree of a Poisson(1) branching process conditioned to live forever,
but these terms will be
defined later, in Section~\ref{5}.  The theorem is stated here only for
the sequence $K_n$, but is in fact true for a more general class of
sequences, which includes $B_n$.

\section{Spanning trees and random walks} \label{2}

Unless $G$ is a very small graph, it is virtually impossible to 
list all of its spanning trees.  For example, if $G = K_n$ is 
the complete graph on $n$ vertices, then the number of spanning 
trees is $n^{n-2}$ according to the well known Pr\"ufer bijection 
\cite{St}.  If $n$ is much bigger than say 20, this is too
many to be enumerated even by the snazziest computer that ever will be.
Luckily, there are shortcuts which enable us to compute probabilities
such as $\P (e \in \Tree)$ without actually enumerating all spanning
trees and counting the proportion containing $e$.  The shortcuts
are based on a close correspondence between spanning trees and
random walks, which is the subject of this section.

\subsection{Simple random walk} \label{2.1}

Begin by defining a simple random walk on $G$.  To avoid obscuring
the issue, we will place extra assumptions on the graph $G$ and later
indicate how to remove these.  In particular, in addition to assuming
that $G$ is finite and connected, we will often suppose that it is $D$-regular 
for some positive integer $D$, which means that every vertex has 
precisely $D$ edges incident to it.  Also suppose that $G$ is simple, 
i.e. it has no self-edges or parallel edges (different edges connecting the
same pair of vertices).  For any vertex $x \in G$, define a simple random 
walk on $G$ starting at $x$, written $SRW_x^G$, intuitively as follows.
Imagine a particle beginning at time 0 at the vertex $x$.  At each
future time $1,2,3,\ldots$, it moves along some edge, always choosing
among the $D$ edges incident to the vertex it is currently at with
equal probability.  When $G$ is not $D$-regular, the definition will
be the same: each of the edges leading away from the current position
will be chosen with probability $1 / \mbox{degree}(v)$.
This defines a sequence of random positions
$SRW_x^G (0), SRW_x^G (1) , SRW_x^G (2) , \ldots$ which is thus
a random function $SRW_x^G$ (or just $SRW$ if $x$ and $G$ may be
understood without ambiguity) from the nonnegative integers to the
vertices of $G$.  Formally, this random function may be defined by
its finite-dimensional marginals which are given by
$\P (SRW_x^G (0) = y_0 , SRW_x^G (1) = y_1 , \ldots , SRW_x^G (k) = y_k) = 
D^{-k}$ if $y_0 = x$ and for all $i = 1 , \ldots , k$ there is an edge 
 from $y_{i-1}$ to $y_i$, and zero otherwise.  
For an illustration of this definition, let $G$
be the following 3-regular simple graph.

\begin{picture}(200,80)(-40,0) 
\put(10,10) {\line(1,0){90}}
\put(10,10) {\line(0,1){60}}
\put(10,70) {\line(1,0){90}}
\put(100,10) {\line(0,1){60}}
\put(40,40) {\line(1,0){30}}
\put(40,40) {\line(-1,1){30}}
\put(40,40) {\line(-1,-1){30}}
\put(70,40) {\line(1,-1){30}}
\put(70,40) {\line(1,1){30}}
\put(6,72){A}
\put(100,72){B}
\put(100,4){C}
\put(6,4){D}
\put(40,42){E}
\put(67,42){F}
\end{picture}
figure~\thepfigure
\label{pfig2.0}
\addtocounter{pfigure}{1} \\
Consider a simple random walk $SRW_A^G$ starting at the vertex $A$.  
The probability of a particular beginning, say $SRW (1) = B$ and $SRW(2)
= F$ is just $(1/3)^2$.  The random position at time 2, $SRW (2)$, is then 
equal to $F$ with probability $2/9$, since each of the two ways, 
ABF and AEF, of getting to $F$ in two steps has probability $1/9$.

Another variant of random walk we will need is the {\em stationary Markov
chain} corresponding to a simple random walk on $G$.  I will preface this
definition with a quick explanation of Markov chains; since I cannot do
justice to this large topic in two paragraphs, the reader is referred to
\cite{IM}, \cite{Fe} or any other favorite introductory probability 
text for further details.  

A (time-homogeneous) {\em Markov chain} on a finite state space $S$ 
is a sequence of random variables $\{ X_i \}$ taking values in $S$, 
indexed by either the integers or the nonnegative integers and having 
the Markov property: there is a set of transition probabilities
$\{ p(x,y) \, : \, x,y \in S \}$ so that the probability of
$X_{i+1}$ being $y$, conditional upon $X_i = x$, is always equal
to $p(x,y)$ regardless of how much more information about the
past you have.  (Formally, this means $\P (X_{i+1} = y \| X_i = x$
and any values of $X_j$ for $j < i )$ is still $p(x,y)$.)
An example of this is $SRW_x^G$, where $S$ is the set of vertices
of $G$ and $p(x,y) = D^{-1}$ if $x \sim y$
and $0$ otherwise (recall that $x \sim y$ means $x$ is a neighbor of $y$).
The values $p(x,y)$ must satisfy $\sum_y p(x,y) =
1$ for every $x$ in order to be legitimate conditional probabilities.
If in addition they satisfy $\sum_x p(x,y) = 1$ for every $y$, the
Markov chain is said to be {\em doubly stochastic}.  It will be useful
later to know that the Markov property is time-reversible, meaning
if $\{ X_i \}$ is a Markov chaing then so is the sequence
$\{\tilde{X}_i = X_{-i} \}$, and there
are backwards transition probabilities $\tilde{p} (x,y)$
for which $\P (X_{i-1} = y \| X_i = x) = \tilde{p} (x,y)$.

If it is possible eventually
to get from every state in $S$ to every other, then there is a
unique {\em stationary distribution} which is a set of
probabilities $\{ \pi (x) \, : \, x \in S \}$ summing to one and
having the property that $\sum_x \pi (x) p(x,y) = \pi (y)$ for all $y$.
Intuitively, this means that if we build a Markov chain with 
transition probabilities $p(x,y)$ and start it by randomizing $X_0$
so that $\P (X_0 = x) = \pi (x)$ then it will also be true that
$\P (X_i = x) = \pi(x)$ for every $i>0$.  A {\em stationary}
Markov chain is one indexed by the integers (as opposed to just the
positive integers),
in which $\P (X_i = x) = \pi (x)$ for some, hence every $i$.  If
a Markov chain is doubly stochastic, it is easy to check that the
uniform distribution $U$ is stationary:
$$\sum_x U (x) p(x,y) = \sum_x |S|^{-1} p(x,y) = |S|^{-1} = U (y).$$
The stationary distribution $\pi$ is unique (assuming every state can 
get to every other) and is hence uniform over all states.

Now we define a stationary simple random walk on $G$ to be a 
stationary Markov chain with state space $V(G)$ and transition
probabilities $p(x,y) = D^{-1}$ if $x \sim y$ and $0$ otherwise.
Intuitively this can be built by choosing $X_0$ at random 
uniformly over $V(G)$, then choosing the $X_i$ for $i > 0$ by
walking randomly from $X_0$ along the edges and choosing the $X_i$ 
for $i < 0$ also by walking randomly from $X_0$, thinking of this
latter walk as going backwards in time.  (For SRW, $p (x,y) = 
p(y,x) = \tilde{p} (x,y)$ so the walk looks the same backwards
as forwards.)

\subsection{The random walk construction of uniform spanning trees} \label{2.2}

Now we are ready for the random walk construction of uniform random
spanning trees.  What we will actually get is a directed spanning
tree, which is a spanning tree together with a choice of vertex
called the root and an orientation on each edge (an arrow pointing 
along the edge in one of the two possible directions) such that
following the arrows always leads to the root.  Of course a directed
spanning tree yields an ordinary spanning tree if you ignore the
arrows and the root.  Here is an algorithm to generate directed trees
 from random walks.

\noindent{\small GROUNDSKEEPER'S ALGORITHM}
\begin{quotation}
Let $G$ be a finite, connected, $D$-regular, simple
graph and let $x$ be any vertex of $G$.  Imagine that we send the
groundskeeper from the local baseball diamond on a walk along the edges
of $G$ starting from $x$; later we will take to be the walk $SRW_x^G$.
She brings with her the wheelbarrow full 
of chalk used to draw in lines.  This groundskeeper is so eager to
choose a spanning tree for $G$ that she wants to chalk a line over
each edge she walks along.  Of course if that edge, along with the
edges she's already chalked, would form a cycle (or is already
chalked), she is not allowed
to chalk it.  In this case she continues walking that edge but
temporarily -- and reluctantly -- shuts off the flow of chalk.
Every time she chalks a new edge she inscribes an arrow pointing
 from the new vertex back to the old.
\end{quotation}

Eventually every vertex is connected to every other by a chalked
path, so no more can be added without forming a cycle and the
chalking is complete.  It is easy to see that the subgraph consisting
of chalked edges is always a single connected component.  The first
time the walk reaches a vertex $y$, the edge just travelled cannot
form a cycle with the other chalked edges.  Conversely, if the
walk moves from $z$ to some $y$ that has been reached before,
then $y$ is connected to $z$ already by some chalked path, 
so adding the edge $zy$ would create a cycle and is not permitted.  Also 
it is clear that following the arrows leads always to vertices that were
visited previously, and hence eventually back to the root.  Furthermore,
every vertex except $x$ has exactly one oriented edge leading out
of it, namely the edge along which the vertex was first reached.

Putting this all together, we have defined a function -- say $\tau$ -- from 
walks on $G$
(infinite sequences of vertices each consecutive pair connected by
an edge) to directed spanning trees of $G$.  Formally $\tau(y_0 , y_1 , y_2 , 
\ldots)$ is the subgraph $H \subseteq G$ such that if $e$ is an oriented
edge from $w$ to $z$ then
$$ e \in H \Leftrightarrow \mbox{ for some } k > 0, y_k = z , y_{k-1}
= w, \mbox{ and there is no } j < k \mbox{ such that } y_j = z .$$
As an example, suppose $SRW_A^G$ in figure 2.1 begins
ABFBCDAE.  Then applying $\tau$ gives the tree with edges 
BA, FB, CB, DC and EA.  

To be completely formal, I should admit that the groundskeeper's
algorithm never stops if there is a vertex that the walk fails
to hit in finite time.  This is not a problem since we are going to
apply $\tau$ to the path of a $SRW$, and this hits every vertex
with probability one.
As hinted earlier, the importance of this construction is the
following equivalence.

\begin{th} \label{rw construction}
Let $G$ be any finite, connected, $D$-regular, simple graph and let
$x$ be any vertex of $G$.  Run a simple random walk $SRW_x^G$ and
let $\Tree$ be the random spanning tree gotten by ignoring the
arrows and root of the random directed spanning tree $\tau(SRW_x^G)$.
Then $\Tree$ has the distribution of a uniform random spanning tree.
\end{th}
 
To prove this it is necessary to consider a stationary simple random
walk on $G$ ($SSRW^G$).  It will be easy to get back to a $SRW_x^G$ because the
sites visited in positive time by a $SSRW^G$ conditioned on being at
$x$ at time zero form a $SRW_x^G$.  Let $T_n$ be the tree 
$\tau(SSRW (n) , SSRW(n+1) , \ldots)$; in other words, $T_n$ is the 
directed tree gotten by applying the groundskeeper's algorithm to the
portion of the stationary simple random walk from time
$n$ onwards.  The first goal is to show that the random collection of 
directed trees $T_n$ forms a time-homogeneous Markov chain
as $n$ ranges over all integers.  

Showing this is pretty straightforward because the transition probabilities
are easy to see.  First note that if $t$ and $u$ are any two directed
trees on disjoint sets of vertices, rooted respectively at $v$ and $w$,
then adding any arrow from $v$ to a vertex in $u$ combines them into a
single tree rooted at $w$.  Now define two operations on directed spanning
trees of $G$ as follows. 
\\[2ex] {\bf Operation} $F(t,x)$: {\it
Start with a directed tree $t$ rooted at $v$.  Choose one of the
the $D$ neighbors of $v$ in $G$, say $x$.  Take away the edge in
$t$ that leads out of $x$, separating $t$ into two trees, rooted at
$v$ and $x$.  Now add an edge from $v$ to $x$, resulting in a single
tree $F(t,x)$.
}
\\[2ex] {\bf Operation} $F^{-1}(t,w)$: {\it
Start with a directed tree $t$ rooted at $x$.  Choose one of the
the $D$ neighbors of $x$ in $G$, say $w$.  Follow the path from
$w$ to $x$ in $t$ and let $v$ be the last vertex on this path
before $x$.   Take away the edge in $t$ that leads out of $v$, 
separating $t$ into two trees, rooted at $x$ and $v$.  Now add an 
edge from $x$ to $w$, resulting in a single directed tree $F^{-1} (t,w)$.
}

It is easy to see that these operations really are inverse to
each other, i.e. if $t$ is rooted at $v$ then $F^{-1} (F(t,x),w) =t$
for any $x \sim v$, where $w$ is the other endpoint of the edge
leading out of $x$ in $t$.  Here is a pictorial example.

\setlength{\unitlength}{1pt}
\begin{picture}(300,140)(-10,0)
\put(100,100){\circle{2}}
\put(95,103){v}
\put(120,80){\circle{2}}
\put(120,80){\vector(-1,1){20}}
\put(80,80){\circle{2}}
\put(80,80){\vector(1,1){20}}
\put(60,60){\circle{2}}
\put(60,60){\vector(1,1){20}}
\put(53,63){w}
\put(40,40){\circle{2}}
\put(40,40){\vector(1,1){20}}
\put(33,43){x}
\put(20,20){\circle{2}}
\put(20,20){\vector(1,1){20}}
\put(0,0){\circle{2}}
\put(0,0){\vector(1,1){20}}
\put(40,0){\circle{2}}
\put(40,0){\vector(-1,1){20}}
\put(350,100){\circle{2}}
\put(345,103){v}
\put(370,80){\circle{2}}
\put(370,80){\vector(-1,1){20}}
\put(330,80){\circle{2}}
\put(330,80){\vector(1,1){20}}
\put(310,60){\circle{2}}
\put(310,60){\vector(1,1){20}}
\put(303,63){w}
\put(290,120){\circle{2}}
\put(283,123){x}
\put(270,100){\circle{2}}
\put(270,100){\vector(1,1){20}}
\put(250,80){\circle{2}}
\put(250,80){\vector(1,1){20}}
\put(290,80){\circle{2}}
\put(290,80){\vector(-1,1){20}}
\put(350,100){\vector(-3,1){60}}
\put(100,40){The tree $t$}
\put(250,40){The tree $F(t,x)$}
\end{picture}  
figure~\thepfigure
\label{pfig2.1}
\addtocounter{pfigure}{1}

I claim that for any directed trees $t$ and $u$,
the backward transition probability $\tilde{p}(t,u)$ is equal to $D^{-1}$
if $u = F(t,x)$ for some $x$ and zero otherwise.  To see this,
it is just a matter of realizing where the operation $F$ comes from.
Remember that $T_n$ is 
just $\tau(SSRW (n) , SSRW(n+1) , \ldots)$, so in particular the root
of $T_n$ is $SSRW (n)$.  Now $SSRW$ really is a Markov chain.
We already know that $\P (SSRW (n-1) = x \| SSRW (n) = v)$ is $D^{-1}$
if $x \sim v$ and zero otherwise.  Also, this is
unaffected by knowledge of $SSRW(j)$ for any $j > n$.  Suppose
it turns out that $SSRW(n-1)=x$.  Then knowing only $T_n$ and
$x$ (but not the values of $SSRW(j)$ for $j > n$) it is possible
to work out what $T_{n-1}$ is.  Remember that $T_n$ and $T_{n-1}$ come
 from applying $\tau$ to the respective sequences $SSRW(n) , SSRW(n+1) ,
\ldots$ and $SSRW(n-1) , SSRW(n), \ldots$ whose only difference
is that the second of these has an extra $x$ tacked on the beginning.
Every time the first sequence reaches a vertex for the first time,
so does the second, unless that vertex happens to be $x$.  So the
$T_{n-1}$ has has all the oriented edges of $T_n$ except the
one out of $x$.  What it has instead is an oriented edge from 
$v$ to $x$, chalked in by the groundskeeper at her very first step.
Adding in the edge from $v$ to some neighbor $x$ and erasing 
the edge out of $x$ yields precisely $F(t,x)$.
So we have shown that $T_{n-1} = F(T_n , SSRW (n-1))$.  But $SSRW(n-1)$
is uniformly distributed among the neighbors of $SSRW(n)$ no matter
what other information we know about the future.  This proves the
claim and the time-homogeneous Markov property.  

The next thing to show is that the stationary distribution is uniform
over all directed trees.  As we've seen, this would follow if
we knew that $\{ T_n \}$ was doubly stochastic.  Since $p(t,u)$
is $D^{-1}$ whenever $u = F(t,x)$ for some $x$ and zero otherwise,
this would be true if for every tree $u$ there are
precisely $D$ trees $t$ for which $F(t,x) = u$ for some $x$.
But the trees $t$ for which $F(t,x) = u$ for some $x$ are
precisely the trees $F^{-1} (u,x)$ for some neighbor $x$ of the
root of $u$, hence there are $D$ such trees and transition
probabilities for $SSRW$ are doubly stochastic.

Now that the stationary distribution for $\{ T_n \}$ has been shown 
to be uniform, the proof of Theorem~\ref{rw construction} is almost
done.  Note that the event $SSRW (0) = x$ is the same 
as the event of $\tau(SSRW (0), SSRW (1) , \ldots)$ being rooted at $x$.
Since $SRW_x^G$ is just $SSRW$ conditioned on $SSRW(0)=x$, 
$T_0 (SRW_x^G)$ is distributed as a uniform directed spanning
tree conditioned on being rooted at $x$.  That is to say, $T_0(SRW_x^G)$
is uniformly distributed over all directed spanning trees rooted
at $x$.  But ordinary spanning trees are in a one to one correspondence
with directed spanning trees rooted at a fixed vertex $x$, the
correspondence being that to get from the ordinary tree to the directed
tree you name $x$ as the root and add arrows that point toward $x$.
Then the tree $\Tree$ gotten from $T_0 (SRW_x^G)$ by ignoring
the root and the arrows is uniformly distributed over all ordinary
spanning trees of $G$, which is what we wanted to prove.    $\Cox$

\subsection{Weighted graphs} \label{2.3}

It is time to remove the extra assumptions that $G$ is $D$-regular
and simple.  It will make sense later to generalize from graphs to
weighted graphs, and since the generalization of Theorem~\ref{rw
construction} is as easy for weighted graphs as for unweighted graphs, 
we may as well introduce weights now.  

A weighted graph is just a graph to each edge $e$ of which is
assigned a positive real number called its weight and written 
$w(e)$.  Edge weights are not allowed to be zero, though one may
conceptually identify a graph with an edge of weight zero with
the same graph minus the edge in question.  An unweighted graph may be
thought of as a graph with all edge weights equal to one, as will be clear
 from the way random trees and random walks generalize.  Write $d(v)$ for the
sum of the weights of all edges incident to $v$.  Corresponding
to the old notion of a uniform random spanning tree is the
{\em weight-selected} random spanning tree ($WST$).  A $WST$,
$\Tree$ is defined to have 
$$ \P (\Tree = t) = { \prod_{e \in t} w(e) \over \sum_u
\prod_{e \in u} w(e) } $$
so that the probability of any individual tree is proportional
to its weight which is by definition the product of the weights of its edges.  

Corresponding to a simple random walk from a vertex $x$ is the
weighted random walk from $x$, $WRW_x^G$ which is a Markov
Chain in which the transition probabilities from a vertex
$v$ are proportional to the weights of the edges incident 
to $v$ (among which the walk must choose).  Thus if
$v$ has two neighbors $w$ and $x$, and there are four edges
incident to $v$ with respective weights $1,2,3$ and $4$ that
connect $v$ respectively to itself, $w$, $x$ and $x$, then
the probabilities of choosing these four edges are respectively
$1/10, 2/10, 3/10$ and $4/10$.  Formally, the probability of walking
along an edge $e$ incident to the current position $v$ is given
by $w(e) / d(v)$.  The bookkeeping is a little
unwieldly since knowing the set of vertices $WRW (0), WRW(1), \ldots$
visited by the $WRW$ does not necessarily determine which edges
were travelled now that the graph is not required to be simple.
Rather than invent some clumsy {\em ad hoc} notation to include
the edges, it is easier just to think that a $WRW$ includes
this information, so it is \underline{not} simply given
by its positions $WRW (j) : j \geq 0$, but that we will refer
to this information in words when necessary.  If $G$ is a connected
weighted graph then $WRW^G$ has a unique stationary distribution
denoted by positive numbers $\pi^G (v)$ summing to one.  This 
will not in general be uniform, but its existence is enough to
guarantee the existence of a stationary Markov chain with the
same transition probabilities.  We call this stationary Markov chain 
$SWRW$ the few times the need arises.  The new and improved
theorem then reads:
\begin{th} \label{weighted rw construction}
Let $G$ be any finite, connected weighted graph and let
$x$ be any vertex of $G$.  Run a weighted random walk $WRW_x^G$ and
let $\Tree$ be the random spanning tree gotten by ignoring the
arrows and root of the random directed spanning tree $\tau(WRW_x^G)$.
Then $\Tree$ has the distribution of $WST$.
\end{th}

The proof of Theorem~\ref{rw construction} serves for
Theorem~\ref{weighted rw construction} with a few alterations.  These
will now be described, thought not much would be lost by taking
these details on faith and skipping to the next section.

The groundskeeper's algorithm is unchanged with
the provision that the $WRW$ brings with it the information of which
edge she should travel if more than one edge connects $WRW(i)$ to
$WRW(i+1)$ for some $i$.  The operation to get from the directed
tree $T_n$ to a candidate for $T_{n-1}$ is basically the same
only instead of there being $D$ choices for how to do this there
is one choice for each edge incident to the root $v$ of $T_n$: choose
such an edge, add it to the tree oriented from $v$ to its other
endpoint $x$ and remove the edge out of $x$.  It is easy to see
again that $\{ T_n \}$ is a time-homogeneous Markov chain with
transition probability from $t$ to $u$ zero unless $u$ can be
gotten from $t$ by the above operation, and if so the probability
is proportional to the weight of the edge that was added in
the operation.  (This is because if $T_n = t$ then $T_{n-1} = u$ 
if and only if $u$ can be gotten from this operation and
$WRW$ travelled along the edge added in this operation between
times $n-1$ and $n$.)  

The uniform distribution on vertices is no longer stationary
for $WRW$ since we no longer have $D$-regularity,
but the distribution $\pi (v) = d(v) / \sum_x d(x)$ is
easily seen to be stationary: start a $WRW$ with $WRW (0)$ having
distribution $\pi$; then
\begin{eqnarray*}
\P (WRW (1) = v) & = & \sum_x \P (WRW (0) = x \mbox{ and }
    WRW(1) = v) \\[2ex]
& = & \sum_x {d(x) \over \sum_y d(y)} \left ( \sum_{e 
    \mbox{\scriptsize ~connecting $x$ to }v} w(e)/d(x) \right ) \\[2ex]
& = & {1 \over \sum_y d(y)} \sum_{e \mbox{\scriptsize ~incident to }v} w(e)
     \\[2ex]
& = & \pi (v) .
\end{eqnarray*}

The stationary distribution $\pi$ for the Markov chain $\{ T_n \}$ 
gives a directed tree $t$ rooted at $v$ probability 
$$\pi (t) = K d(v) \prod_{e \in t} w(e), $$
where $K = (\sum_t d(\mbox{root}(t)) 
\prod_{e \in t} w(e))^{-1}$ is a normalizing constant.   
If $t$, rooted at $v$, can go to $u$, rooted at $x$, by adding an
edge $e$ and removing the edge $f$, then $\pi (u) / \pi (t)
= d(x) w(e) / d(v) w(f)$.  To verify that $\pi$ is a stationary
distribution for $T_n$ write $\C (u)$ for the class of trees
 from which it is possible to get to $u$ in one step and for 
each $t \in \C (u)$ write $v_t, e_t$ and $f_t$ for the root of
$t$, edge added to $t$ to get $u$ and edge taken away from
$t$ to get $u$ respectively.  If $u$ is rooted at $x$, then
\begin{eqnarray*}
\P (T_{n-1} = u) & = & \sum_t \P (T_n = t \mbox{ and }
    T_{n-1} = u) \\[2ex]
& = & \sum_{t \in \C (u)} \pi (t) w(e_t) / d(v_t) \\[2ex]
& = & \sum_{t \in \C (u)} [\pi (u) d(v_t) w(f_t) /d(x) w(e_t)] w(e_t)
    / d(v_t) \\[2ex]
& = & \pi (u) \sum_{t \in \C (u)} w(f_t) / d(x) \\[2ex]
& = & \pi (u) ,
\end{eqnarray*}
since as $t$ ranges over all trees that can get to $u$, $f_t$
ranges over all edges incident to $x$.  

Finally, we have again that $\tau(WRW_x^G (0))$ is distributed
as $\tau(SWRW^G (0))$ conditioned on having root $x$, and since
the unconditioned $\pi$ is proportional to $d(x)$ times
the weight of the tree (product of the edge weights), the factor of $d$ 
is constant and $\P (\tau(WRW_x^G(0)) = t)$ is proportional to 
$\prod_{e \in t} w(e)$ for any $t$ rooted at $x$.  Thus 
$\tau(WRW_x^G(0))$ is distributed identically to $WST$.    $\Cox$

\subsection{Applying the random walk construction to our model} \label{2.4}

Although the benefit is not yet clear, we have succeeded in translating
the question of determining $\P (e \in \Tree)$ from a question
about uniform spanning trees to a question about simple random walks.
To see how this works, suppose that $e$ connects the vertices 
$x$ and $y$ and generate a uniform spanning tree by the random
walk construction starting at $x$: $\Tree =$ the tree gotten from
$\tau(SRW_x^G)$ by ignoring the root and arrows.  If $e \in \Tree$
then its orientation in $\tau(SRW_x^G)$ must be from $y$ to $x$,
and so $e \in \Tree$ if and only if $SRW (k-1) = x$ where
$k$ is the least $k$ for which $SRW (k) = y$.  In other words,
\begin{equation} \label{rw prob 1}
\P (e \in \Tree) = \P (\mbox{first visit of } SRW_x^G \mbox{ to }
y \mbox{ is along } e) .
\end{equation}
The computation of this random walk probability turns out to be tractable.

More important is the fact that this may be iterated to get
probabilities such as $\P (e,f \in \Tree)$.  This requires two
more definitions.  If $G$ is a finite connected graph and $e$
is an edge of $G$ whose removal does not disconnect $G$, then
the {\em deletion} of $G$ by $e$ is the graph $G \setminus e$
with the same vertex set and the same edges minus $e$.  If $e$
is any edge that connects distinct vertices $x$ and $y$,
then the {\em contraction} of $G$ by $e$ is the graph $G/e$ 
whose vertices are the vertices of $G$ with $x$ and $y$ replaced
by a single vertex $x*y$.  There is an edge $\rho (f)$ of $G/e$ for every
edge of $f$ of $G$, where if one or both endpoints of $f$
is $x$ or $y$ then that endpoint is replaced by $x*y$ in $\rho (f)$.
We write $\rho (z)$ for the vertex corresponding to $z$ in this
correspondence, so $\rho (x) = \rho (y) = x*y$ and $\rho (z) = z$
for every $z \neq x,y$.  The following example shows $G_1$ and
$G_1 / e_4$.  The edge $e_4$ itself maps to a self-edge under $\rho$,
$e_5$ becomes parallel to $e_6$ and $D$ and $E$ map to $D*E$. \\
\setlength{\unitlength}{2pt}
\begin{picture}(200,80)(60,0)
\put(100,50){\circle{2}}
\put(98,53){A}
\put(140,50){\circle{2}}
\put(140,53){B}
\put(140,10){\circle{2}}
\put(140,3){C}
\put(100,10){\circle{2}}
\put(98,3){D}
\put(70,30){\circle{2}}
\put(62,28){E}
\put(100,50){\line(1,0){40}}
\put(120,53){$e_1$}
\put(100,50){\line(0,-1){40}}
\put(102,30){$e_6$}
\put(140,50){\line(0,-1){40}}
\put(142,30){$e_2$}
\put(140,10){\line(-1,0){40}}
\put(120,3){$e_3$}
\put(100,50){\line(-3,-2){30}}
\put(77,42){$e_5$}
\put(100,10){\line(-3,2){30}}
\put(77,12){$e_4$}
\put(155,30){$G_1$}
\put(200,50){\circle{2}}
\put(198,53){A}
\put(240,50){\circle{2}}
\put(240,53){B}
\put(240,10){\circle{2}}
\put(240,3){C}
\put(200,10){\circle{2}}
\put(196,4){\scriptsize D*E}
\put(200,50){\line(1,0){40}}
\put(220,53){$e_1$}
\put(204,30){$e_6$}
\put(240,50){\line(0,-1){40}}
\put(242,30){$e_2$}
\put(240,10){\line(-1,0){40}}
\put(220,3){$e_3$}
\put(200,30){\oval(5,40)}
\put(185,30){$e_5$}
\put(200,2){\circle {16}}
\put(185,3){$e_4$}
\put(255,30){$G_1 / e_4$}
\end{picture} 
figure~\thepfigure
\label{pfig2.2}
\addtocounter{pfigure}{1}

It is easy to see that successive deletions and contractions may be
performed in any order with the same result.  If $e_1 , \ldots , e_r$
are edges of $G$ whose joint removal does not disconnect $G$ then the
successive deletion of these edges is permissible.  Similarly 
if $\{ f_1, \ldots , f_s \}$ is a set of edges of $G$ that contains
no cycle, these edges 
may be successively contracted and the graph $G \setminus e_1, \ldots , e_r /
f_1 , \ldots , f_s$ is well-defined.  It is obvious that the spanning 
trees of $G \setminus e$ are just those spanning trees of $G$ that do not 
contain $e$.  Almost as obvious is a one to one correspondence between 
spanning trees of $G$ containing $e$ and spanning trees of $G/e$: 
if $t$ is a spanning tree of $G$ containing $e$ then there is a 
spanning tree of $G/e$ consisting of $\{ \rho (f) : f \neq e \in t \}$.

To translate $\P (e,f \in \Tree)$ to the random walk setting, 
write this as $\P (e \in \Tree) \P (f \in \Tree \| e \in \Tree)$.  
The first term has already been translated.  The conditional
distribution of a uniform random spanning tree given that it contains
$e$ is just uniform among those trees containing $e$, which
is just $\P_{G/e} (\rho (f) \in \Tree)$ where the subscript $G/e$
refers to the fact that $\Tree$ is now taken to be a uniform
random spanning tree of $G/e$.  If $f$ connects $z$ and $x$
then this is in turn equal to $\P (SRW_{\rho (x)}^{G/e}
\mbox{ first hits } \rho (z) \mbox{ along } \rho(f))$.  Both the terms
have thus been translated; in general it should be clear how this may
be iterated to translate the probability of any elementary 
event, $\P (e_1 , \ldots , e_r \in \Tree \mbox{ and } f_1 , \ldots ,
f_s \notin \Tree)$ into a product of random walk probabilities. 
It remains to be seen how these probabilities may be calculated.

\section{Random walks and electrical networks} \label{3}

Sections~\ref{3.1} - \ref{3.3} contain a development of 
the connection between random walks and electrical networks.
The right place to read about this is in \cite{DS}; what
you will see here is necessarily a bit rushed.  Sections~\ref{3.5}
and~\ref{3.6} contain similarly condensed material from
other sources.

\subsection{Resistor circuits} \label{3.1}

The electrical networks we discuss will have only two kinds of 
elements: resistors and voltage sources.  Picture the resistors as
straight pieces of wire.  A resistor network will be built by soldering
resistors together at their endpoints.  That means that a diagram
of a resistor network will just look like a finite graph with
each edge bearing a number: the resistance.
Associated with every resistor network $H$ is a weighted
graph $G_H$ which looks exactly like the graph just mentioned except
that the weight of an edge is not the resistance but the
{\em conductance}, which is the reciprocal of the resistance.
The distinction between $H$ and $G_H$ is only necessary while we are
discussing precise definitions and will then be dropped.
A voltage source may be a single battery that provides a specified
voltage difference (explained below) across a specified pair
of vertices or is may be a more complicated device to hold various
voltages fixed at various vertices of the network.
Here is an example of a resistor network on a familiar graph,
with a one volt battery drawn as a dashed box.  Resistances on 
the edges (made up arbitrarily) are given in ohms.

\setlength{\unitlength}{2pt}
\begin{picture}(180,80)
\put(100,50){\circle{2}}
\put(90,53){A}
\put(140,50){\circle{2}}
\put(142,53){B}
\put(140,10){\circle{2}}
\put(140,3){C}
\put(100,10){\circle{2}}
\put(98,2){D}
\put(70,30){\circle{2}}
\put(62,28){E}
\put(100,50){\line(1,0){40}}
\put(120,44){1}
\put(100,50){\line(0,-1){40}}
\put(102,30){2}
\put(140,50){\line(0,-1){40}}
\put(142,30){2}
\put(140,10){\line(-1,0){40}}
\put(120,1){.5}
\put(100,50){\line(-3,-2){30}}
\put(77,42){1}
\put(100,10){\line(-3,2){30}}
\put(77,12){1}
\put (106,56) {\dashbox{2}(28,14){1 V}}
\put (106,63) {\line (-1,-2){6}}
\put (100,72) {\scriptsize +}
\put (134,63) {\line (1,-2){6}}
\put (137,72) {-}
\end{picture}
figure~\thepfigure
\label{pfig3.1}
\addtocounter{pfigure}{1}

The electrical properties of such a network are given by Kirchoff's
laws.  For the sake of exposition I will give the laws numbers,
although these do not correspond to the way Kirchoff actually 
stated the laws.  The first law is that every vertex of the network has
a voltage which is a real number.  
The second law gives every oriented edge (resistor) a current.
Each edge has two possible orientations.  Say an edge connects
$x$ and $y$.  Then the current through the edge is a real number
whose sign depends on which orientation you choose for the edge.
In other words, the current $I(\xy)$ that flows from $x$ to $y$
is some real number and the current $I(\yx)$ is its negative.
(Note though that the weights $w(e)$
are always taken to positive; weights are functions of unoriented
edges, whereas currents are functions of oriented edges.)
If $I (e)$ denotes the current along an oriented edge $e = \xy$, 
$V(x)$ denotes the voltage at $x$ and $R(e)$ denotes the resistance 
of $e$, then quantatively, the second law says
\begin{equation} \label{KL1}
I(\xy) = [V(x) - V(y)] R(e)^{-1} . 
\end{equation}
Kirchoff's third law is that the total current flowing into
a vertex equals the total current flowing out, or in other
words
\begin{equation} \label{KL1.5}
\sum_{y \sim x} I(\xy) = 0 .
\end{equation}
This may be rewritten using~(\ref{KL1}).
Recalling that in the weighted graph $G_H$, the weight $w(e)$
is just $R(e)^{-1}$ and that $d(v)$ denotes the sum of $w(e)$
over edges incident to $v$, we get at every vertex $x$ an equation
\begin{equation} \label{KL2}
0 = \sum_{y \sim x} [V(x) - V(y)] w(xy) =
   V(x) d(x) - \sum_{y \sim x} V(y) w(xy) .
\end{equation}

Since a voltage source may provide current, this may fail to hold at
any vertex connected to a voltage source.  The above laws are
sufficient to specify the voltages of the network -- and hence the
currents -- except that a constant may be added to all the voltages
(in other words, it is the voltage
differences that are determined, not the absolute voltages).  
In the above example the voltage difference across
$AB$ is required to be one.  Setting the voltage at $B$ to
zero (since the voltages are determined only up to an additive constant)
the reader may check that the voltages at $A,C,D$ and $E$ are
respectively $1,4/7,5/7$ and $6/7$ and the currents through
$AB,AE,ED,AD,DC,CB$ are respectively $1,1/7,1/7,1/7,2/7,2/7$.

\subsection{Harmonic functions} \label{3.2}

The voltages in a weighted graph $G$ (which we are now identifying with
the resistor network it represents) under application of a voltage
source are calculated by finding a solution to Kirchoff's laws on $G$
with specified boundary conditions.  For each vertex $x$ there is an unknown
voltage $V(x)$.  There is also a linear equation for every 
vertex not connected to a voltage source, and an equation given by
the nature of each voltage source.  Will these always be 
enough information so that Kirchoff's laws have a unique solution?
The answer is yes and it is most easily seen in the context of harmonic
functions.  \footnote{There is also the question of whether any
solution exists, but addressing that would take us too far afield.
If you aren't convinced of its existence on physical grounds, wait
until the next subsection where a probabilistic interpretation for the 
voltage is given, and then deduce existence of a solution from 
the fact that these probabilities obey Kirchoff's laws.}

If $f$ is a function on the vertices of a weighted graph $G$, define
the {\em excess} of $f$ at a vertex $v$, written $\lap f(v)$ by
$$ \lap f (v) = \sum_{y \sim v} [f(v) - f(y)] w(vy) .$$
You can think of $\lap$ as an operator that maps functions $f$ to
other functions $\lap f$ that is a discrete analog of
the Laplacian operator.   
A function $f$ from the vertices of a finite weighted graph $G$
to the reals is said to be {\em harmonic} at a vertex
$v$ if and only if $\lap f(v) = 0$.
Note that for any function $f$, the sum of the excesses
$\sum_{v \in G} \lap f(v) = 0$, since each $[f(x)-f(y)]w(xy)$
cancels a $[f(y)-f(x)]w(yx)$ due to $w(xy) = w(yx)$.
To see what harmonic functions are intuitively, consider the special case 
where $G$ is unweighted, i.e. all of the edge weights are one.  Then a
function is harmonic if and only if its value at a vertex $x$ is the
average of the values at the neighbors of $x$.  In the weighted
case the same is true, but with a weighted average!  Here is an easy
but important lemma about harmonic functions.

\begin{lem}[Maximum principle] \label{max principle}
Let $V$ be a function on the vertices of a finite connected weighted
graph, harmonic everywhere except possibly at vertices of
some set $X = \{ x_1 , \ldots , x_k \}$.
Then $V$ attains its maximum and minimum on $X$.  If
$V$ is harmonic everywhere then it is constant.
\end{lem}

\noindent{Proof:}  Let $S$ be the set of vertices where $V$ attains
its maximum.  Certainly $S$ is nonempty.  If $x \in S$ has a neighbor
$y \notin S$ then $V$ cannot be harmonic at $x$ since $V(x)$ would
then be a weighted average of values less than or equal to $V(x)$
with at least one strictly less.  In the case where $V$ is harmonic
everywhere, this shows that no vertex in $S$ has a neighbor not in
$S$, hence since the graph is connected every vertex is in $S$ and
$V$ is constant.  Otherwise, suppose $V$ attains its maximum at
some $y \notin X$ and pick a path connecting
$y$ to some $x \in X$.  The entire path must then be in $S$ up until 
and including the first vertex along the path at which $V$ is not 
harmonic.  This is some $x' \in X$.  The argument for the minimum is just
the same.   $\Cox$

Kirchoff's third law~(\ref{KL2}) says that the voltage function is 
harmonic at every $x$ not connected to a voltage source.  
Suppose we have a voltage source that provides a fixed voltage
at some specified set of vertices.  Say for concreteness that
the vertices are $x_1, \ldots , x_k$ and the voltages produced
at these vertices are $c_1 , \ldots , c_k$.  We now show that
Kirchoff's laws determine the voltages everywhere else, i.e. there
is at most one solution to them.

\begin{th} \label{unique fix V}
Let $V$ and $W$ be real-valued functions on the vertices of a finite 
weighted graph $G$.  Suppose that $V(x_i) = W(x_i) = c_i$ for some set
of vertices $x_1 , \ldots , x_k$ and $1 \leq i \leq k$ and that
$V$ and $W$ are harmonic at every vertex other than $x_1 , \ldots , x_k$.
Then $V=W$.
\end{th}

\noindent{Proof:}  Consider the function $V-W$.  It is easy to check
that being harmonic at $x$ is a linear property, so $V-W$ is
harmonic at every vertex at which both $V$ and $W$ are harmonic.
Then by the Maximum Principle, $V-W$ attains its maximum and minimum
at some $x_i$.  But $V-W = 0$ at every $x_i$, so $V-W \equiv 0$.
$\Cox$

Suppose that instead of fixing the voltages at a number of points,
the voltage source acts as a current source and supplies a 
fixed amount of current $I_i$ to vertices $x_i , 1 \leq i \leq k$.
This is physically reasonable only if $\sum_{i=1}^k I_i = 0$.
Then a net current of $I_i$ will have to flow out of each $x_i$
into the network.  Using~(\ref{KL1}) gives 
$$I_i = \sum_{y \sim x} w(x,y) (V(x) - V(y)) =
\lap V(x) .$$
 From this it is apparent that the assumption $\sum_i I_i = 0$
is algebraically as well as physically necessary since the excesses must 
sum to zero.  Kirchoff's laws also determine the voltages (up to an 
additive constant) of a network with current sources, as we now show.

\begin{th} \label{unique fix I}
Let $V$ and $W$ be real-valued functions on the vertices of a finite 
weighted graph $G$.  Suppose that $V$ and $W$ both have excess
$c_i$ at $x_i$ for some set of vertices $x_i$ and reals $c_i$ , $1 \leq
i \leq k$.  Suppose also that $V$ and $W$ are harmonic elsewhere.  Then
$V=W$ up to an additive constant.
\end{th}

\noindent{Proof:}  Excess is linear, so the excess of $V-W$ is the
excess of $V$ minus the excess of $W$.  This is zero everywhere, so
$V-W$ is harmonic everywhere.  By the Maximum Principle, $V-W$
is constant.    $\Cox$

\subsection{Harmonic random walk probabilities} \label{3.3}

Getting back to the problem of random walks, suppose $G$ is a finite
connected graph and $x,a,b$ are vertices of $G$.  Let's say that I want
to calculate the probability that $SRW_x$ reaches $a$ before $b$.
Call this probability $h_{ab} (x)$.  It is not immediately obvious 
what this probability is, but we can get an equation by watching
where the random walk takes its first step.  Say the neighbors of $x$
are $y_1 , \ldots , y_d$.  Then $\P (SRW_x (1) = y_i) = d^{-1}$ for
each $i \leq d$.  If we condition on $\P (SRW_x (1) = y_i)$ then
the probability of the walk reaching $a$ before $b$ is (by the
Markov property) the same as if it had started out at $y_i$.  This
is just $h_{ab} (y_i)$.  Thus
\begin{eqnarray*}
h_{ab} (x) & = & \sum_i \P (SRW_x (1) = y_i) h_{ab} (y_i) \\[2ex]
&=& d^{-1} \sum_i h_{ab} (y_i) .
\end{eqnarray*}
In other words, $h_{ab}$ is harmonic at $x$.  Be careful though, if
$x$ is equal to $a$ or $b$, it doesn't make sense to look one step 
ahead since $SRW_x (0)$ already determines whether the walk hit
$a$ or $b$ first.  In particular, $h_{ab} (a) = 1$ and $h_{ab} (b) = 0$,
with $h_{ab}$ being harmonic at every $x \neq a,b$.  

Theorem~\ref{unique fix V} tells us that there is only one such
function $h_{ab}$.  This same function solves Kirchoff's laws for
the unweighted graph $G$ with voltages at $a$ and $b$ fixed at
$1$ and $0$ respectively.  In other words, the probability of $SRW_x$
reaching $a$ before $b$ is just the voltage at $x$ when a one volt
battery is connected to $a$ and $b$ and the voltage at $b$ is 
taken to be zero.  If $G$ is a weighted graph, we can use a similar
argument: it is easy to check that the first-step transition 
probabilities $p(x,y) = w(\xy) / \sum_z w(\xz)$ show that
$h_{ab} (x)$ is harmonic in the sense of weighted graphs.  Summarizing
this:
\begin{th} \label{elec rw}
Let $G$ be a finite connected weighted graph.  Let $a$ and $b$ be 
vertices of $G$.  For any vertex $x$, the probability of $SRW_x^{G}$
reaching $a$ before $b$ is equal to the voltage at $x$ in $G$
when the voltages at $a$ and $b$ are fixed at one and zero volts
respectively.
\end{th}
Although more generality will not be needed we remark that this same
theorem holds when $a$ and $b$ are taken to be sets of vertices.
The probability of $SRW_x$ reaching a vertex in $a$ before reaching
a vertex in $b$ is harmonic at vertices not in $a \cup b$, is zero
on $b$ and one on $a$.  The voltage when vertices in $b$ are held at
zero volts and vertices in $a$ are held at one volt also satisfies
this, so the voltages and the probabilities must coincide.

Having given an interpretation of voltage in probabilistic terms,
the next thing to find is a probabilistic interpretation of the
current.  The arguments are similar so they will be treated briefly; a
more detailed treatment appears in \cite{DS}.  First we will need
to find an electrical analogue for the numbers $u_{ab} (x)$ which
are defined probabilistically as the expected number of times 
a $SRW_a$ hits $x$ before the first time it hits $b$.  This is
defined to be zero for $x=b$.  For any $x \neq a,b$, let 
$y_1, \ldots , y_r$ be the neighbors of $x$.  Then the number
of visits to $x$ before hitting $b$ is the sum over $i$ of the number
of times $SRW_a$ hits $y_i$ before $b$ and goes to $x$ on the
next move (the walk had to be somewhere the move before it hit $x$).
By the Markov property, this quantity is $u_{ab} (y_i) p(y_i,x) 
= u_{ab} (y_i) w(\xy_i) / d(y_i)$.  Letting $\phi_{ab} (z)$
denote $u_{ab} (z) / d(z)$ for any $z$, this yields
$$\phi_{ab} (x) = d(x) u_{ab} (x) = \sum_i u_{ab} (y_i) w(\xy_i) / d(y_i)
   = \sum_i w(\xy_i) \phi_{ab} (y_i) .$$
In other words $\phi_{ab}$ is harmonic at every $x \neq a,b$. 
Writing $K_{ab}$ for $\phi_{ab} (a)$ we then have that
$\phi_{ab}$ is $K_{ab}$ at $a$, zero at $b$ and harmonic elsewhere,
hence it is the same function as the the voltage induced by 
a battery of $K_{ab}$ volts connected to $a$ and $b$, with the
voltage at $b$ taken to be zero.  Without yet knowing what $K_{ab}$
is, this determines $\phi_{ab}$ up to a constant multiple.  This
in turn determines $u_{ab}$, since $u_{ab} (x) = d(x) \phi_{ab} (x)$.

Now imagine that we watch $SRW_a$ to see when it crosses
over a particular edge $\xy$ and count plus one every
time it crosses from $x$ to $y$ and minus one every time it crosses
 from $y$ to $x$.  Stop counting as soon as the walk hits $b$.
Let $H_{ab}(\xy)$ denote the expected number of
signed crossings.  ($H$ now stands for harmonic, not for the name of a
resistor network.)  We can calculate $H$ in terms of $u_{ab}$
by counting the plusses and the minuses separately.  The expected
number of plus crossings is just the expected number of times
the walk hits $x$, mulitplied by the probability on each of
these occasions that the walk
crosses to $y$ on the next move.  This is $u_{ab} (x) 
w(\xy) / d(x)$.  Similarly the expected number of minus
crossings is $u_{ab} (y) w(\xy) / d(y).$  Thus
\begin{eqnarray*}
H_{ab} (\xy) & = & u_{ab} (x) w(\xy) / d(x) - u_{ab} (y) w(\xy) / 
   d(y) \\[2ex]
& = & w(\xy) [ \phi_{ab} (x) - \phi_{ab} (y) ] .
\end{eqnarray*}
But $\phi_{ab} (x) - \phi_{ab} (y)$ is just the voltage difference
across $\xy$ induced by a $K_{ab}$-volt battery across $a$ and $b$.
Using~(\ref{KL1}) and $w(\xy) = R(\xy)^{-1}$ shows that the expected
number of signed crossings of $\xy$ is just the current induced 
in $\xy$ by a $K_{ab}$-volt battery connected to $a$ and $b$.  
A moment's thought shows that the expected number of signed
crossings of all edges leading out of $a$ must be one, since the
walk is guaranteed to leave $a$ one more time than it returns to
$a$.  So the current supplied by the $K_{ab}$-volt battery must
be one amp.  Another way of saying this is that 
\begin{equation} \label{eq2}
\lap \phi_{ab} = \dd_a - \dd_b . 
\end{equation}
Instead of worrying about what $K_{ab}$ is, we may
just as well say that the expected number of crossings of
$\xy$ by $SRW_a$ before hitting $b$ is the current induced
when one amp is supplied to $a$ and drawn out at $b$.

\subsection{Electricity applied to random walks applied to spanning trees}
\label{3.4}

Finally we can address the random walk question that relates to
spanning trees.  In particular, the claim that the probability
in equation~(\ref{rw prob 1}) is tractable will be borne out several
different ways.  First we will see how the probability may be
``calculated'' by an analog computing device, namely a resistor
network.  In the next subsection, the computation will be carried out 
algebraically and very neatly, but only for particularly nice 
symmetric graphs.  At the end of the section, a universal method 
will be given for the computation which is a little messier.  Finally in 
Section~\ref{4} the question of the individual probabilities in~(\ref{rw prob
1}) will be avoided altogether and we will see instead how values
for these probabilities (wherever they might come from) determine
the probabilities for all contractions and deletions of the graph
and therefore determine all the joint probabilities $\P (e_1 , \ldots
, e_k \in \Tree)$ and hence the entire measure.  

Let $e = \xy$ be any edge of a finite connected weighted graph $G$.
Run $SRW_x^G$ until it hits $y$.  At this point either the walk just
moved along $e$ from $x$ to $y$ -- necessarily for the first time -- and
$e$ will be in the tree $\Tree$ given by $\tau(SRW_x^G)$, or else
the walk arrived at $y$ via a different edge in which case the walk
never crossed $e$ at all and $e \notin \Tree$.  In either case the
walk never crossed from $y$ to $x$ since it stops if it hits $y$.
Then the expected number of signed crossings of $e=\xy$ by $SRW_x$ 
up to the first time it hits $y$ is equal to the probability of 
first reaching $y$ along $e$ which equals $\P (e \in \Tree)$.
Putting this together with the electrical interpretation of signed
crossings give
\begin{th} \label{current}
$\P (e \in \Tree) = $ the fraction of the current that goes through
edge $e$ when a battery is hooked up to the two endpoints of $e$.
\end{th}
$\Cox$

This characterization leads to a proof of Theorem~\ref{neg cor}
provided we are willing to accept a proposition that is 
physically obvious but not so easy to prove, namely
\begin{th}[Rayleigh's monotonicity law] \label{Rayleigh}
The effective resistance of a circuit cannot increase when
a new resistor is added.
\end{th}
The reason this is physically obvious is that adding a new 
resistor provides a new path for current
to take while allowing the current still to flow through 
all the old paths.  Theorem~\ref{neg cor} says that the conditional
probability of $e \in \Tree$ given $f \in \Tree$ must be less
than or equal to the unconditional probability.  Using 
Theorem~\ref{current} and the fact that the probabilities 
conditioned on $f \notin \Tree$ are just the probabilities for 
$WST$ on $G \setminus f$, this boils down to showing that the fraction 
of current flowing directly across $e$ is no greater on $G$
than it is on $G \setminus f$.  The battery across $e$ meets two parallel
resistances: $e$ and the effective resistance of the rest of $G$.
The fraction of current flowing through $e$ is inversely proportional
to the ratio of these two resistances.  Rayleigh's theorem says
that the effective resistance of the rest of $G$ including
$f$ is at most the effective resistance of $G \setminus f$, so the fraction
flowing through $e$ on $G$ is at most the fraction flowing
through $e$ on $G \setminus f$.
In Section~\ref{4}, a proof will be given that does not rely
on Rayleigh.

\subsection{Algebraic calculations for the square lattice} \label{3.5}

If $G$ is a finite graph, then the functions from the vertices of $G$
to the reals form a finite-dimensional real vector space.  The operator
$\lap$ that maps a function $V$ to its excess is a linear operator
on this vector space.  In this
language, the voltages in a resistor network with one 
unit of current supplied to $a$ and drawn out at $b$ are the unique
(up to additive constant) function $V$ that solves $\lap V = 
\dd_a - \dd_b$.  Here $\dd_x$ is the function that is one at $x$ 
and zero elsewhere.  This means that $V$ can be calculated simply
by inverting $\lap$ in the basis $\{ \dd_x ; x \in G \}$.  Although
$\lap$ is technically not invertible, its nullspace has dimension one
so it can be inverted on a set of codimension one.  A classical
determination of $V$ for arbitrary graphs is carried out in the
next subsection.  The point of this subsection is to show how the 
inverse can be obtained in a simpler way for nice graphs.  

The most general ``nice'' graphs to which the method will apply are
the infinite $\Z^d$-periodic lattices.  Since in this article I am
restricting attention to finite graphs, I will not attempt to be general
but will instead show a single example.  The reader may look in
\cite{BP} for further generality.  The example considered here is
the square lattice.  This is just the graph you see on a piece of
graph paper, with vertices at each pair of integer coordinates and four
edges connecting each point to its nearest neighbors.  The exposition
will be easiest if we consider a finite square piece of this and impose
wrap-around boundary conditions.  Formally, let $T_n$ (T for torus) be
the graph whose vertices are pairs of integers $\{ (i,j) : 0 
\leq i,j \leq n-1 \}$ and for which two points are connected if and only
if they agree in one component and differ by one mod $n$ in the other
component.  Here is a picture of this with $n=3$ and the broken edges
denoting edges that wrap around to the other side of the graph.  The
graph is unweighted (all edge weights are one.)

\begin{picture}(150,90) 
\put (20,20) {\circle {2}} 
\put (20,20) {\line (1,0) {20}} 
\put (20,20) {\line (0,1) {20}} 
\put (20,40) {\circle {2}}
\put (20,40) {\line (1,0) {20}} 
\put (20,40) {\line (0,1) {20}} 
\put (20,60) {\circle {2}}
\put (20,60) {\line (1,0) {20}} 
\put (20,60) {\dashbox(0,15)}
\put (40,20) {\circle {2}}
\put (40,20) {\line (1,0) {20}} 
\put (40,20) {\line (0,1) {20}} 
\put (40,40) {\circle {2}}
\put (40,40) {\line (1,0) {20}} 
\put (40,40) {\line (0,1) {20}} 
\put (40,60) {\circle {2}}
\put (40,60) {\line (1,0) {20}} 
\put (40,60) {\dashbox{1}(0,15)}
\put (60,20) {\circle {2}}
\put (60,20) {\dashbox{1}(15,0)}
\put (60,20) {\line (0,1) {20}}
\put (60,40) {\circle {2}}
\put (60,40) {\dashbox{1}(15,0)}
\put (60,40) {\line (0,1) {20}}
\put (60,60) {\circle {2}}
\put (60,60) {\dashbox{1}(0,15)}
\put (60,60) {\dashbox{1}(15,0)}
\end{picture}
figure~\thepfigure
\label{pfig3.2}
\addtocounter{pfigure}{1}

Let $\zeta = e^{2\pi i/n}$ denote the first $n^{th}$ root of unity.
To invert $\lap$ we exhibit its eigenvectors.  Since the vector
space is a space of functions, the eigenvectors are called
eigenfunctions.  For each pair of integers $0 \leq k,l \leq n-1$ 
let $f_{kl}$ be the function on the vertices of $T_n$ defined by 
$$f_{kl} (i,j) = \zeta^{ki+lj} .$$
If you have studied group representations, you will recognize $f_{kl}$
as the representations of the group $T_n = (\Z / n\Z)^2$ and in fact 
the rest of this section may be restated more compactly in terms of
characters of this abelian group.

It is easy to calculate 
\begin{eqnarray*}
\lap f_{kl} (i,j) & = & 4 \zeta^{ki+lj} - \zeta^{ki+l(j+1)} - \zeta^{ki+l(j-1)} 
   - \zeta^{k(i+1)+lj} - \zeta^{k(i-1)+lj} \\[2ex]
& = & \zeta^{ki+lj} (4 - \zeta^k - \zeta^{-k} - \zeta^l - \zeta^{-l} )
   \\[2ex]
& = & \zeta^{ki+lj} (4 - 2 \cos (2\pi k/n) - 2 \cos (2 \pi l/n)) .
\end{eqnarray*}
Since the multiplicative factor $(4 - 2 \cos (2\pi k/n) - 2 \cos 
(2 \pi l/n))$ does not depend on $i$ or $j$, this shows that
$f_{kl}$ is indeed an eigenfunction for $\lap$
with eigenvalue $\lambda_{kl} = 4 - 2 \cos (2\pi k/n) - 2 \cos 
(2 \pi l/n)$.

Now if $\{ v_k \}$ are eigenvectors for some linear operator
$A$ with eigenvalues $\{ \lambda_k \}$, then for any constants 
$\{ c_k \}$, 
\begin{equation} \label{invert eig}
A^{-1} (\sum_k c_k v_k) = \sum_k \lambda_k^{-1} c_k v_k .
\end{equation}
If some $\lambda_k$ is equal to zero, then the range of $A$ does not
include vectors $w$ with $c_k \neq 0$, so $A^{-1} w$ does not exist for
such $w$ and indeed the formula blows up due to the $\lambda_k^{-1}$.
In our case $\lambda_{kl} =  4 - 2 \cos (2\pi k/n) - 2 \cos 
(2 \pi l/n) = 0$ only when $k=l=0$.  
Thus to calculate $\lap^{-1} (\dd_a - \dd_b)$ we need to 
figure out coefficents $c_{kl}$ for which $\dd_a - \dd_b = 
\sum_{kl} c_{kl} f_{kl}$ and verify that $c_{00} = 0$.
For this puropose, it is fortunate that 
the eigenfunctions $\{ f_{kl} \}$ are actually a unitary basis
in the inner product $<f,g> = \sum_{ij} f(i,j) \overline{g(i,j)}$.
You can check this by writing 
$$< f_{kl} , f_{k'l'} > = \sum_{ij} \zeta^{ki+lj} \overline
{\zeta^{k'i+l'j}};$$
elementary algebra show this to be one if $k=k'$ and $l=l'$ and zero 
otherwise, which what it means to be unitary.  Unitary bases are
great for calculation because the coefficients $\{c_{kl} \}$ of 
any $V$ in a unitary eigenbasis $\{ f_{kl} \}$ are 
given by $c_{kl} = <V , f_{kl}>$.  In our case, this means
$c_{kl} = \sum_{ij} V(i,j) \overline{f_{kl} (i,j)}$.  Letting 
$a$ be the vertex $(0,0)$, $b$ be the vertex $(1,0)$ and $V = \dd_a 
- \dd_b$, this gives $c_{kl} = 1 - {\overline{\zeta}}^k$ and hence
$$ \dd_a - \dd_b = \sum_{k,l} (1 - {\overline{\zeta}}^k ) f_{kl} .$$

We can now plug this into equation~(\ref{invert eig}), since
clearly $c_{00} = 0$.  This gives
\begin{eqnarray}
\lap V & = & \dd_a - \dd_b \nonumber \\[2ex]
&&\nonumber \\
\Leftrightarrow \;\;\; V(i,j) & = & cf_{00}(i,j) + \sum_{(k,l) \neq
    (0,0)} (1 - \zeta^k) \lambda_{kl}^{-1} f_{kl}(i,j) \nonumber \\[2ex]
& = & c + \sum_{(k,l) \neq (0,0)} { 1 - \zeta^k \over 4 - 
    2 \cos (2 \pi k / n) - 2 \cos (2 \pi l/n)} \zeta^{ki+lj} .  \label{eq1}
\end{eqnarray}
This sum is easy to compute exactly and to approximate efficiently
when $n$ is large.  In particular as $n \rightarrow \infty$ the
sum may be replaced by an integral which by a small miracle admits
an exact computation.  Details of this may be found in \cite[page 148]{Sp}.
You may check your arithmetic against mine by using~(\ref{eq1}) to derive the
voltages for a one volt battery placed across the bottom left
edge $e$ of $T_3$ and across the bottom left edge $e'$ of $T_4$:
$$\begin{array}{ccc} 5/8 & 3/8 & 1/2 \\ 5/8 & 3/8 & 1/2 \\ 1 & 0 & 1/2
  \end{array} \hspace{2in} \begin{array}{cccc} 56/90 & 34/90 & 40/90 & 50/90 \\ 
  50/90 & 40/90 & 42/90 & 48/90 \\ 56/90 & 34/90 & 40/90 & 50/90 \\
  1 & 0 & 34/90 & 56/90 \end{array} \hfill $$

Section~\ref{5} shows how to put these numbers to good use, but
we can already make one calculation based on Theorem~\ref{current}.
The four currents flowing out of the bottom left vertex under 
the voltages shown are given by the voltage differences: $1, 3/8, 1/2$
and $3/8$.  The fraction of the current flowing directly through
the bottom left edge $e$ is $8/18$, and according to Theorem~\ref{current},
this is $\P (e \in \Tree)$.  An easy way to see this is right is by
the symmetry of the graph $T_3$.  Each of the 18 edges should
be equally likely to be in $\Tree$, and since every spanning tree
has 8 edges, the probability of any given edge being in the tree
must be $8/18$.

\subsection{Electrical networks and spanning trees} \label{3.6}

The order in which topics have been presented so far makes sense
 from an expository viewpoint but is historically backwards.  The
first interest in enumerating spanning trees came from problems in 
electrical network theory.  To set the record straight and also to close
the circle of ideas 
\begin{quote}
spanning trees $\rightarrow$ random walks $\rightarrow$ electrical 
networks $\rightarrow$ spanning trees 
\end{quote}
I will spend a couple of paragraphs on this remaining connection.

Let $G$ be a finite weighted graph.  Assume there are no voltage
sources and the quantity of interest is the {\em effective resistance}
between two vertices $a$ and $b$.  This is defined to be the voltage
it is necessary to place across $a$ and $b$ to induce a unit current 
flow.  A classical theorem known to Kirchoff is:
\begin{th} \label{eff resistance}
Say $s$ is an $a,b$-spanning bitree if $s$ is a spanning forest
with two components, one containing $a$ and the other containing $b$.
The effective resistance between $a$ and $b$ may be computed from
the weighted graph $G$ by taking the quotient $N / D$ where
$$D = \sum_{\mbox{\scriptsize spanning trees }t} \hspace{1in} 
    \left ( \prod_{e \in t} w(e) \right )$$
is the sum of the weights of all spanning trees of $G$ and
$$N = \sum_{a,b\mbox{\scriptsize -spanning bitrees }s} \hspace{1in} 
    \left ( \prod_{e \in s} w(e) \right )$$
is the analogous sum over $a,b$-spanning bitrees.    $\Cox$
\end{th}

To see that how this is implied by Theorem~\ref{current} and
equation~(\ref{rw prob 1}), imagine adding an extra one ohm resistor
 from $a$ to $b$.  The probability of this edge being chosen in
a $WST$ on the new graph is by definition given by summing the 
weights of trees containing the new edge and
dividing by the total sum of the weights of all spanning trees.  
Clearly $D$ is the
sum of the weights of trees not containing the extra edge.  But
the trees containing the extra edge are in one-to-one correspondence
with $a,b$-spanning bitrees (the correspondence being to remove
the extra edge).  The extra edge has weight one, so the sum of the
weights of trees that do contain the extra edge is $N$ and the probability
of a $WST$ containing the extra edge is $N/(N+D)$.  By equation~(\ref{rw
prob 1}) and
Theorem~\ref{current}, this must then be the fraction of current flowing
directly through the extra edge when a battery is placed across $a$ and
$b$.  Thinking of the new circuit as consisting of the extra edge in
parallel with $G$, the fractions of the current passing through the
two components are proportional to the inverses of their resistances,
so the ratio of the resistance of the extra edge to the rest of
the circuit must be $D:N$.  Since the extra edge has resistance one,
the effective resistance of the rest of the circuit is $N/D$.

The next problem of course was to efficiently evaluate the sum of 
the weights of all spanning trees of a weighted graph.  The solution
to this problem is almost as well known and can be found, among
other places in \cite{CDS}.
\begin{th}[Matrix-Tree Theorem] \label{matrix tree}
Let $G$ be a finite, simple, connected, weighted graph 
and define a matrix indexed by the vertices of $G$ by letting 
$M(x,x) = d(x)$, $M(x,y) = - w(\xy)$ if $x$ and $y$ are connected
by an edge, and $M(x,y) = 0$ otherwise.  Then for any vertex $x$, 
the sum of the weights of all spanning trees of $G$ is equal to 
the determinant of the matrix gotten from $M$ by deleting by the
row and column corresponding to $x$.  
\end{th}
The matrix $M$ is nothing but a representation of $\lap$ with
respect to the basis $\{ \dd_x \}$.  Recalling that the problem
essentially boils down to inverting $\lap$, the only other
ingredient in this theorem is the trick of inverting the
action of a singular matrix on an element on its range by inverting
the largest invertible principal minor of the matrix.  Details can
be found in \cite{CDS}.   $\Cox$

\section{Transfer-impedances} \label{4}

In the last section we saw how to calculate $\P (e \in \Tree)$ in 
several ways: by Theorems~\ref{current} or~\ref{eff resistance} in 
general and by equations such as~(\ref{eq1}) in particularly symmetric
cases.  By repeating the calculations in Theorem~\ref{current} 
and~\ref{eff resistance} for contractions and deletions of a graph
(see Section~\ref{2.4}),
we could then find enough conditional probabilities to determine the
probability of any elementary event $\P (e_1 , \ldots , e_r \in \Tree
\mbox{ and } f_1 , \ldots , f_s \notin \Tree)$.  Not only is this
inefficient, but it fails to apply to the symmetric case of
equation~(\ref{eq1}) since contracting or deleting the graph breaks
the symmetry.  The task at hand is to alleviate this problem by showing
how the data we already know how to get -- current flows on $G$ --
determine the current flows on contractions and deletions of $G$ and 
thereby determine all the elementary probabilities for $WST$ on $G$.
This will culminate in a proof of Theorem~\ref{exist matrix}, which
encapsulates all of the necessary computation into a single determinant.

\subsection{An electrical argument} \label{4.1}

To keep notation to a minimum this subsection will only deal with 
unweighted, $D$-regular graphs.  
Begin by stating explicitly the data that will be used to determine all
other probabilities.  
For oriented edges $e = \xy$ and $f = \zw$ in a finite connected
graph $G$, define the {\em transfer-impedance} $H (e,f) =
\phi_{xy} (z) - \phi_{xy} (w)$ which is equal to
the voltage difference across $f$, $V(z) - V(w)$, when one
amp of current is supplied to $x$ and drawn out at $y$.  
We will assume knowledge of
$H (e,f)$ for every pair of edges in $G$ (presumably via some
analog calculation, or in a symmetric case by equation~(\ref{eq1})
or something similar) and show how to derive all other 
probabilities from these transfer-impedances.

Note first that $H(e,e)$ is the voltage across $e$ for a unit current
flow supplied to one end of $e$ and drawn out of the other.
This is equal to the current flowing directly
along $e$ under a unit current flow and is thus $\P (e \in \Tree)$.
The next step is to try a computation involving a single contraction.
For notation, recall the map $\rho$ which projects vertices and edges of
$G$ to vertices and edges of $G/f$.
Fix edges $e=\xy$ and $f=\zw$ and let $\{ V(v) : v \in G/f\}$ 
be the voltages we need to solve for: voltages at vertices of $G/f$ 
when a unit current is supplied to $\rho (x)$ and drawn out at 
$\rho (y)$.  As we have seen, this means $\lap V (v) = +1, -1$ or $0$
according to whether $v = x,y$ or neither.
Suppose we lift this to a function $\VV$ on the
vertices of $G$ by letting $\VV (x) = V( \rho (x))$.  Let's calculate
the excess $\lap \VV$ of $\VV$.  Each edge of $G$ corresponds to
an edge in $G/f$, so for any $v \neq z,w$ in $G$, $\lap \VV (v)
= \lap V(\rho (v))$; this is equal to $+1$ if $v=x$, $-1$ if $v=y$
and zero otherwise.  Since $\rho$ maps both $z$ and $w$ onto
the same vertex $v*w$, we can't tell what the $\lap \VV$ is at 
$z$ or $w$ individually, but $\lap \VV (z) + \lap \VV (w)$ will 
equal $\lap V (z*w)$ which will equal $+1$ if $z$ or $w$ coincides
with $x$, $-1$ if $z$ or $w$ coincides with $y$ and zero otherwise
(or if both coincide!).  The last piece of information we have 
is that $\VV (z) = \VV (w)$.  Summarizing,
\begin{quote}
$(i)$ $\lap \VV = \dd_x - \dd_y + c (\dd_z - \dd_w)$; \\

$(ii)$ $\VV (z) = \VV (w)$ ,
\end{quote}
where $c$ is some unknown constant.  To see that this uniquely defines
$\VV$ up to an additive constant, note that the difference between
any two such functions has excess $c (\dd_z - \dd_w)$ for some $c$,
hence by the maximum principle reaches its maximum and minimum on 
$\{ z , w \}$; on the other hand the values at $z$ and $w$ are equal,
so the difference is constant.  

Now it is easy to find $\VV$.  Recall from equation~(\ref{eq2}) that 
$\phi$ satisfies $\lap \phi_{ab} = \dd_a - \dd_b$.  The function
$\VV$ we are looking for is then $\phi_{xy} + c \phi_{zw}$ where $c$ is
chosen so that 
$$\phi_{xy} (z) + c \phi_{zw} (z) = \phi_{xy} (w) + c \phi_{zw} (w) . $$
In words, $\VV$ gives the voltages for a battery supplying unit
current in at $x$ and out at $y$ plus another battery across $z$ and $w$
just strong enough to equalize the voltages at $z$ and $w$.  How
strong is that?  The battery supplying unit current to $x$ and $y$
induces by definition a voltage $H(\xy , \zw)$ across $z$ and $w$. 
To counteract that, we need a $-H(\xy , \zw)$-volt battery across
$z$ and $w$.  Since supplying one unit of current in at $z$ and out at $w$
produces a voltage across $z$ and $w$ of $H(\zw , \zw)$, the current 
supplied by the counterbattery must be $c = -H(\xy , \zw) / H(\zw , \zw)$.
We do not need to worry about $H (\zw , \zw)$ being zero since this
means that $\P (f \in \Tree) = 0$ so we shouldn't be conditioning
on $f \in \Tree$.  Going back to the original problem, 
\begin{eqnarray*}
\P (e \in \Tree \| f \in \Tree) & = & V (\rho (x)) - V (\rho(y)) \\[2ex]
& = & \VV (x) - \VV (y) \\[2ex]
& = & H (\xy , \xy) + H(\zw , \xy) { - H (\xy , \zw ) \over H(\zw , \zw)}
    \\[2ex]
& = & {H (\xy , \xy) H(\zw , \zw) - H(\xy , \zw) H(\zw , \xy) \over
    H(\zw , \zw)} .
\end{eqnarray*}
Multiplying this conditional probability by the unconditional
probability $\P (f \in \Tree)$ gives the probability of both $e$ and $f$
being in $\Tree$ which may be written as
$$\P (e,f \in \Tree) = \left | \begin{array}{cc} H(\xy,\xy) &
   H(\xy , \zw) \\ H(\zw , \xy) & H(\zw , \zw) \end{array} \right | .$$
Thus $\P (e,f \in \Tree) = \det M(e,f)$ where $M$ is the matrix 
of values of $H$ as in Theorem~\ref{exist matrix}.  

Theorem~\ref{exist matrix} has in fact now been proved for $r = 1,2$.
The procedure for general $r$ will be similar.  Write $\P (e_1 , 
\ldots e_r \in \Tree)$ as a product of conditional probabilities
$\P (e_i \in \Tree \| e_{i+1} , \ldots , e_r \in \Tree)$.  Then
evaluate this conditional probability by solving for voltages on
$G/e_{i+1} \cdots e_r$.  This is done by placing batteries
across $e_1 , \ldots , e_r$ so as to equalize voltages across
all $e_{i+1} , \ldots , e_r$ simultaneously.  Although in the
$r=2$ case it was not necessary to worry about dividing by 
zero, this problem does come up in the general case which
causes an extra step in the proof.  We will now summarily generalize
the above discussion on how to solve for voltages on contractions 
of a graph and then forget about electricity altogether.

\begin{lem} \label{gen contraction}
Let $G$ be a finite $D$-regular connected graph and let $f_1 , \ldots ,
f_r$ and $e = \xy$ be edges of $G$ that form no cycle.  Let $\rho$ be the 
map from
$G$ to $G/f_1 \ldots f_r$ that maps edges to corresponding edges and
maps vertices of $G$ to their equivalence classes under the relation
of being connected by edges in $\{ f_1 , \ldots , f_r \}$.  Let
$\VV$ be a function on the vertices of $G$ such that 
\begin{quote}
$(i)$ If $\zw = f_i$ for some $i$ then $\VV (z) = \VV (w)$ ; \\

$(ii)$ $\sum_{z \in \rho^{-1} (v)} \lap \VV (z) = +1$ if $\rho (x)
=v$, $-1$ if $\rho (y) = v$ and zero otherwise.
\end{quote}
If $\Tree$ is a uniform spanning tree for $G$ then $\P (e \in 
\Tree \| f_1 , \ldots , f_r \in \Tree ) = \VV (x) - \VV (y)$.
\end{lem}

\noindent{Proof:}  As before, we know that $\P (e \in \Tree \| f_1 ,
\ldots , f_r \in \Tree)$ is given by $V(\rho (x)) - V(\rho (y))$
where $V$ is the voltage function on $G/f_1 \cdots f_r$ for
a unit current supplied in at $x$ and out at $y$.  Defining $\VV (v)$
to be $V(\rho (v))$, the lemma will be proved if we can show that
$\VV$ is the unique function on the vertices of $G$ 
satisfying~$(i)$ and~$(ii)$.  Seeing that $\VV$ satisfies~$(i)$ 
and~$(ii)$ is the same as before.  Since $\rho$ provides
a one to one correspondence between edges of $G$ and edges of
$G/f_1 , \ldots , f_r$, the excess of $\VV$ at vertices of
$\rho^{-1}(v)$ is the sum over edges leading out of vertices
in $\rho^{-1} (v)$ of the difference of $\VV$ across that edge,
which is the sum over edges leading out of $\rho (v)$ of the
difference of $V$ across that edge; this is the excess of $V$ at
$\rho (v)$ which is $=1,-1$ or $0$ according to whether $x$ or
$y$ or neither is in $\rho^{-1} (v)$.

Uniqueness is also easy.  If $\WW$ is
any function satisfying $(i)$, define a function
$W$ on the vertices of $G/f_1 \cdots f_r$ by $W(\rho (v)) = \WW (v)$.
If $\WW$ satisfies $(ii)$ as well then it is easy to check that
$W$ satisfies $\lap W = \dd_{\rho (x)} - \dd_{\rho(y)}$ so that 
$W = V$ and $\WW = \VV$.   $\Cox$

\subsection{Proof of the transfer-impedance theorem} \label{4.2}

First of all, though is is true that the function $H$ in the
previous subsection and the statement of the theorem is
symmetric, I'm not going to include a proof -- nothing else
we talk about relies on symmetry of $H$ and a proof may be found
in any standard treatment of the Green's function, such as 
\cite{Sp}. Secondly, it is easiest to reduce the problem to the
case of $D$-regular graphs immediately so as to be able to use the
previous lemma.  Suppose $G$ is any finite connected graph.  Let $D$ be 
the maximum degree of any vertex in $G$ and to any vertex of
lesser degree $k$, add $D-k$ self-edges.  The resulting
graph is $D$-regular (though not simple) and furthermore it has
the same spanning trees as $G$.  To prove Theorem~\ref{exist matrix}
for finite connected graphs, it therefore suffices to prove
the theorem for finite, connected, $D$-regular graphs.  Restating
what is to be proved:
\begin{th} \label{transfer thm}
Let $G$ be any finite, connected, $D$-regular graph and let
$\Tree$ be a uniform random spanning tree of $G$.  Let
$H(\xy , \zw)$ be the voltage induced across $\zw$ when one
amp is supplied from $x$ to $y$.  Then for any $e_1 , \ldots , e_r
\in G$, 
$$ \P (e_1 , \ldots , e_r \in \Tree) = \det M(e_1 , \ldots , e_r)$$
where $M(e_1 , \ldots , e_r)$ is the $r$ by $r$ matrix whose
$i,j$-entry is $H(e_i , e_j)$.  
\end{th}

The proof is by induction on $r$.  We have already proved it
for $r = 1,2$, so now we assume it for $r-1$ and try to prove it
for $r$.  There are two cases.  The first possibility is that
$\P (e_1 , \ldots , e_{r-1} \in \Tree) = 0$.  This means that
no spanning tree of $G$ contains $e_1 , \ldots , e_r$ which means
that these edges contain some cycle.  Say the cycle is $e_{n(0)},
\ldots , e_{n(k-1)}$ where there are vertices $v(i)$ for which
$e_{n(i)}$ connects $v(i)$ to $v(i+1 \mbox{ mod } k)$.  For any
vertices $x,y$, $\phi_{xy}$ is the unique solution up
to an additive constant of $\lap \phi_{xy} = \dd_x - \dd_y$.  Thus
$\lap \left ( \sum_{i=0}^{k-1} \phi_{v(i)\,v(i+1 \mbox{ mod } k)} \right ) = 0$
which means that $\sum_{i=0}^{k-1} \phi_{v(i) \, v(i+1 \mbox{ mod } k)}$
is constant.  Then for any $\xy$,  
\begin{eqnarray*}
&&\sum_{i=0}^{k-1} H(e_{n(i)} , \xy) \\[2ex]
& = & \sum_{i=0}^{k-1} \phi_{v(i) \, v(i+1 \mbox{ mod } k)} (x) -
    \sum_{i=0}^{k-1} \phi_{v(i) \, v(i+1 \mbox{ mod } k)} (y) \\[2ex]
& = & 0 .
\end{eqnarray*}
This says that in the matrix $M(e_1 , \ldots , e_r)$, the rows
$n(1), \ldots , n(k)$ are linearly dependent, summing to zero.
Then $\det M(e_1 , \ldots , e_r) = 0$ which is certainly the probability
of $e_1 , \ldots , e_r \in \Tree$.

The second possibility is that $\P (e_1 , \ldots , e_{r-1} \in \Tree) 
\neq 0$.  We can then write 
\begin{eqnarray*}
&& \P (e_1 , \ldots , e_r \in \Tree) \\[2ex]
& = & \P (e_1 , \ldots e_{r-1} \in \Tree) \P (e_r \in \Tree \| e_1 , 
   \ldots e_{r-1} \in \Tree) \\[2ex]
& = & \det M(e_1 , \ldots , e_{r-1}) \P (e_r \in \Tree \| e_1 , 
   \ldots e_{r-1} \in \Tree)
\end{eqnarray*}
by the induction hypothesis.  To evaluate the last term we look for a
function $\VV$ satisfying the conditions of Lemma~\ref{gen contraction}
with $e_r$ instead of $e$ and  $e_1 , \ldots, e_{r-1}$ instead of $f_1 
, \ldots , f_r$.
For $i \leq r-1$, let $x_i$ and $y_i$ denote the vertices connected 
by $e_i$.  For any $v \in G/e_1 \cdots e_{r-1}$ and any $i \leq r-1$,
$\sum_{z \in \rho^{-1} (v)} \lap \phi_{x_i y_i} (z) =
\sum_{z \in \rho^{-1} (v)} \dd_{x_i} (z) - \dd_{y_i} (z)$
which is zero since the class $\rho^{-1} (v)$ contains both
$x_i$ and $y_i$ or else contains neither.  The excess of $\phi_{x_r
y_r}$ summed over $\rho^{-1} (v)$ is just $1$ if $\rho (x_r) =
v$, $-1$ if $\rho (y_r) = v$ and zero otherwise.  By linearity
of excess, this implies that the sum of $\phi_{x_r y_r}$ with
any linear combination of $\{ \phi_{x_i y_i} : i \leq r-1 \}$ 
satisfies $(ii)$ of the lemma.  

Satisfying part $(i)$ is then a matter of choosing the right linear
combination, but the lovely thing is that we don't have to actually
compute it!  We do need to know it exists and here's the argument
for that.  The $i^{th}$ row of $M(e_1 , \ldots , e_r)$ lists the
values of $\phi_{x_i y_i} (x_j) - \phi_{x_i y_i} (y_j)$ as $j$
runs from 1 to $r$.  Looking for $c_1 , \ldots , c_{r-1}$ 
such that $\phi_{x_r y_r} + \sum_{i=1}^{r-1} \phi_{x_i y_i}$
is the same on $x_j$ as on $y_j$ for $j \leq r-1$ is the
same as looking for $c_i$ for which the $r^{th}$ row of 
$M$ plus the sum of $C_i$ times the $i^{th}$ row of $M$ has
zeros for every entry except the $r^{th}$.  In other words
we want to row-reduce, using the first $r-1$ rows to clear
$r-1$ zeros in the last row.  There is a unique way to do this
precisely when the determinant of the upper $r-1$ by $r-1$
submatrix is nonzero, which is what we have assumed.  So these
$c_1 , \ldots , c_{r-1}$ exist and $\VV (v) = \phi_{x_r y_r} (v)
+ \sum_{i=1}^{r-1} \phi_{x_i y_i} (v)$.  

The lemma tells us that $\P (e_r \in \Tree \| e_1 , \ldots , e_{r-1}
\in \Tree)$ is $\VV (x_r) - \VV (y_r)$.  This is just the $r,r$-entry
of the row-reduced matrix.  Now calculate the determinant of
the row-reduced matrix in two ways.  Firstly, since row-reduction
does not change the determinant of a matrix, the determinant must
still be $\det M(e_1 , \ldots , e_r)$.  On the other hand, since
the last row is all zeros except the last entry, expanding along
the last row gives that the determinant is the $r,r$-entry times
the determinant of the upper $r-1$ by $r-1$ submatrix, which
is just $\P (e_r \in \Tree \| e_1 , \ldots , e_{r-1} \in \Tree)
\det M(e_1 , \ldots , e_{r-1})$.  Setting these two equal gives
$$\P (e_r \in \Tree \| e_1 , \ldots , e_{r-1} \in \Tree)
= \det M(e_1 , \ldots , e_r) / \det M(e_1 , \ldots , e_{r-1}) .$$
The induction hypothesis says that
$$ \P (e_1 , \ldots , e_{r-1} \in \Tree) = \det M(e_1 , \ldots e_{r-1})$$
and multiplying the conditional and unconditional probabilities proves
the theorem.    $\Cox$

\subsection{A few computational examples} \label{4.3}

It's time to take a break from theorem-proving to see how well the
machinery we've built actually works.  A good place to test it is
the graph $T_3$, since the calculations have essentially been
done, and since even $T_3$ is large enough to prohibit enumeration
of the spanning trees directly by hand (you can use the Matrix-Tree
Theorem with all weights one to check that there are 11664 of them).
Say we want to know the probability that the middle vertex $A$
is connected to $B,C$ and $D$ in a uniform random spanning tree
$\Tree$ of $T_3$.\\
\setlength{\unitlength}{2pt}
\begin{picture}(150,90) 
\put (20,20) {\circle {2}} 
\put (20,20) {\line (1,0) {20}} 
\put (20,20) {\line (0,1) {20}} 
\put (20,40) {\circle {2}}
\put(22,42){B}
\put (20,40) {\line (1,0) {20}} 
\put (20,40) {\line (0,1) {20}} 
\put (20,60) {\circle {2}}
\put (20,60) {\line (1,0) {20}} 
\put (20,60) {\dashbox(0,15)}
\put (40,20) {\circle {2}}
\put(42,22){E}
\put (40,20) {\line (1,0) {20}} 
\put (40,20) {\line (0,1) {20}} 
\put (40,40) {\circle {2}}
\put(42,42){A}
\put (40,40) {\line (1,0) {20}} 
\put (40,40) {\line (0,1) {20}} 
\put (40,60) {\circle {2}}
\put (40,60) {\line (1,0) {20}} 
\put(42,62){C}
\put (40,60) {\dashbox{1}(0,15)}
\put (60,20) {\circle {2}}
\put (60,20) {\dashbox{1}(15,0)}
\put (60,20) {\line (0,1) {20}}
\put (60,40) {\circle {2}}
\put(62,42){D}
\put (60,40) {\dashbox{1}(15,0)}
\put (60,40) {\line (0,1) {20}}
\put (60,60) {\circle {2}}
\put (60,60) {\dashbox{1}(0,15)}
\put (60,60) {\dashbox{1}(15,0)}
\end{picture}
figure~\thepfigure
\label{pfig4.1}
\addtocounter{pfigure}{1} \\

We need then to calculate the transfer-impedance matrix for the
edges $AB,AC$ and $AD$.  Let's say we orient them all toward $A$.
The symmetry of $T_3$ under translation and $90^\circ$ rotation 
allows us to rely completely on the voltages calculated at the end 
of~3.5.  Sliding the picture upwards one square and multiplying
the given voltages by $4/9$ to produce a unit current flow from $B$
to $A$ gives voltages
$$\begin{array}{ccc} 5/18 & 3/18 & 4/18 \\ 8/18 & 0 & 4/18 \\
       5/18 & 3/18 & 4/18 \end{array} $$
which gives transfer-impedances $H(BA,BA) = 8/18$, $H(BA,CA) = 3/18$
and $H(BA,DA)= 4/18$.  The rest of the values follow by symmetry,
giving 
$$M(BA,CA,DA) = {1 \over 18} \left ( \begin{array}{ccc} 8 & 3 & 4 \\
   3 & 8 & 3 \\ 4 & 3 & 8 \end{array} \right ) .$$
Applying Theorem~\ref{transfer thm} gives
$\P (BA,CA,DA \in \Tree) = \det M(BA,CA,DA) = {\displaystyle {312 \over 
5832}}$, or in other words just 624 of the 11664 spanning trees of $T_3$ 
contain all these edges.  Compare this to using the Matrix-Tree Theorem
to calculate the same probability.  That does not require the
preliminary calculation of the voltages, but it does require an
eight by eight determinant.  

Suppose we want now to calculate the probability that $A$ is 
a {\em leaf} of $\Tree$, that is to say there is only one edge in 
$\Tree$ incident to $A$.  By symmetry this edge will be $AB$
$1/4$ of the time, so we need to calculate $\P (BA \in \Tree \mbox{ and } 
CA,DA,EA \notin \Tree)$ and then multiply by four.  As 
remarked earlier, we can use inclusion-exclusion to get the
answer.  This would entail writing 
\begin{eqnarray*}
& & \P (BA \in \Tree \mbox{ and } CA,DA,EA \notin \Tree) \\[2ex]
& = & \P (BA \in \Tree) - \P (BA,CA \in \Tree) - \P (BA,DA \in \Tree)
   - \P (BA,EA \in \Tree) \\
&& + \P (BA,CA,DA \in \Tree) + \P (BA,CA,EA \in \Tree) + \P (BA,DA,EA
   \in \Tree)  \\
&& - \P (BA,CA,DA,EA \in \Tree) .
\end{eqnarray*}
This is barely manageable for four edges, and gets exponentially messier
as we want to know about probabilities involving more edges.  Here is
an easy but useful theorem telling how to calculate the probability of a 
general {\em cylinder} event, namely the event that $e_1 , \ldots , e_r$ are 
in the tree, while $f_1 , \ldots , f_s$ are not in the tree.  

\begin{th} \label{incl-excl}
Let $M (e_1 , \ldots , e_k)$ be an $k$ by $k$
transfer-impedance matrix.  Let $M^{(r)}$ be the matrix for which 
$M^{(r)}(i,j) = M(i,j)$ if $i \leq r$ and $M^{(r)}(i,j) = 1 - M(i,j)$ if
$r+1 \leq i \leq k$.  Then $\P (e_1 , \ldots , e_r \in \Tree
\mbox{ and } e_{r+1} , \ldots , e_k \notin \Tree) = \det M^{(r)}$.
\end{th}

\noindent{Proof:}  The proof is by induction on $k-r$.  The initial
step is when $r=k$; then $M^{(r)} = M$ so the theorem
reduces to Theorem~\ref{transfer thm}.  Now suppose the theorem
to be true for $k-r = s$ and let $k-r = s+1$.  Write 
\begin{eqnarray*}
&&\P (e_1 , \ldots , e_r \in \Tree \mbox{ and } e_{r+1} , \ldots , 
e_k \notin \Tree) \\[2ex]
& = & \P (e_1 , \ldots , e_r \in \Tree \mbox{ and } e_{r+2} , \ldots , 
    e_k \notin \Tree) \\
& & -  \P (e_1 , \ldots , e_{r+1} \in \Tree \mbox{ and } e_{r+2} , \ldots , 
    e_k \notin \Tree) \\[2ex]
& = & \det M(e_1 , \ldots , e_r , e_{r+2} , \ldots e_k) -
    \det M(e_1 , \ldots , e_{r+1} , e_{r+2} , \ldots e_k) ,
\end{eqnarray*}
since the induction hypothesis applies to both of the last two
probabilities.  Call these last two matrices $M_1$ and $M_2$.
The trick now is to stick an extra row and column into $M_1$: let $M'$
be $M(e_1 , \ldots , e+k)$ with the $r+1^{st}$ row replaced by
zeros except for a one in the $r+1^{st}$ position.  Then $M'$
is $M_1$ with an extra row and column inserted.  Expanding 
along the extra row gives $\det M' = \det M_1$.  But $M'$ and 
$M_2$ differ only in the $r+1^{st}$ row, so by multilinearity
of the determinant,
$$\det M_1 - \det M_2 = \det M' - det M_2 = \det M''$$
where $M''$ agrees with $M'$ and $M_2$ except that the 
$r+1^{st}$ row is the difference of the $r+1^{st}$ rows of 
$M'$ and $M_2$.  The induction is done as soon as you realize
that $M''$ is just $M^{(r)}$.   $\Cox$

Applying this to the probability of $A$ being a leaf of $T_3$, 
we write 
\begin{eqnarray*}
&&\P (BA \in \Tree \mbox{ and } CA,DA,EA \notin \Tree) \\[2ex]
& = & \det M^{(3)} (BA,CA,DA,EA) \\[3ex]
& = & \left | \begin{array}{cccc} 8/18 & 3/18 & 4/18 & 3/18 \\
   -3/18 & 10/18 & -3/18 & -4/18 \\ -4/18 & -3/18 & 10/18 & -3/18 \\
   -3/18 & -4/18 & -3/18 & 10/18 \end{array} \right | \;\; = \;\; 
    {10584 \over 18^4} \;\; = \;\; {1176 \over 11664}
\end{eqnarray*}
so $A$ is a leaf of $4 \cdot 1176 = 4704$ of the 11664 spanning trees
of $T_3$.  This time, the Matrix-Tree Theorem would have required
evaluation of several different eight by eight determinants.  If
$T_3$ were replaced by $T_n$, the transfer-impedance
calculation would not be significantly harder, but the
Matrix-Tree Theorem would require several $n^2$ by $n^2$ determinants.
If $n$ goes to $\infty$, as it might when calculating some sort
of limit behavior, these large determinants would not be tractable.

\section{Poisson limits} \label{5}

As mentioned in the introduction, the random degree of a vertex
in a uniform spanning tree of $G$ converges in distribution to
one plus a Poisson(1) random
variable as $G$ gets larger and more highly connected.  This section
investigates some such limits, beginning with an example symmetric
enough to compute explicitly.  The reason for this limit may seem
clearer at the end of the section when we discuss a stronger limit
theorem.  Proofs in this section are mostly
sketched since the details occupy many pages in \cite{BP}.

\subsection{The degree of a vertex in $K_n$} \label{5.1}

The simplest situation in which to look for a Poisson limit is
on the complete graph $K_n$.  This is pictured here for $n=8$. \\
\setlength{\unitlength}{2pt}
\begin{picture}(190,120)(-70,-10)
\put(0,24){\circle* {6}}
\put(0,24){\line(0,1){36}}
\put(0,24){\line(2,5){24}}
\put(0,60){\circle* {6}}
\put(0,60){\line(1,1){24}}
\put(0,60){\line(5,2){60}}
\put(24,84){\circle* {6}}
\put(24,84){\line(1,0){36}}
\put(24,84){\line(5,-2){60}}
\put(60,84){\circle* {6}}
\put(60,84){\line(1,-1){24}}
\put(60,84){\line(2,-5){24}}
\put(84,60){\circle* {6}}
\put(84,60){\line(0,-1){36}}
\put(84,60){\line(-2,-5){24}}
\put(84,24){\circle* {6}}
\put(84,24){\line(-1,-1){24}}
\put(84,24){\line(-5,-2){60}}
\put(60,0){\circle* {6}}
\put(60,0){\line(-1,0){36}}
\put(60,0){\line(-5,2){60}}
\put(24,0){\circle* {6}}
\put(24,0){\line(-1,1){24}}
\put(24,0){\line(-2,5){24}}
\put(0,24){\line(1,0){84}}
\put(0,60){\line(1,0){84}}
\put(24,0){\line(0,1){84}}
\put(60,0){\line(0,1){84}}
\put(0,24){\line(1,1){60}}
\put(24,0){\line(1,1){60}}
\put(60,0){\line(-1,1){60}}
\put(84,24){\line(-1,1){60}}
\put(-3,24){\line(5,2){90}}
\put(-3,60){\line(5,-2){90}}
\put(24,-3){\line(2,5){36}}
\put(60,-3){\line(-2,5){36}}
\end{picture}
figure~\thepfigure
\label{pfig5.1}
\addtocounter{pfigure}{1}

Calculating the voltages for a complete graph is particularly easy
because of all the symmetry.  Say the vertices of $K_n$ are called
 $v_1 , \ldots , v_n$, and put a one volt battery across $v_1$ and
$v_2$, so $V(v_1) = 1$ and $V(v_2) = 0$.  By Theorem~\ref{elec rw}, 
the voltage at any other vertex $v_j$ is equal to the probability
that $SRW_{v_j}^{K_n}$ hits $v_1$ before $v_2$.  This is clearly equal
to $1/2$.  The total current flow out of $v_1$ with these voltages
is $n/2$, since one amp flows along the edge to $v_2$ and $1/2$
amp flows along each of the $n-2$ other edges out of $v_1$.  Multiplying
by $2/n$ to get a unit current flow gives voltages
$$V(v_i) = \left \{ \begin{array}{ccl} ~~2/n~~ & : & ~~i = 1 \\
   ~~0~~ & : & ~~i = 2 \\ ~~1/n~~ & & \mbox{  otherwise}.
   \end{array} \right.  $$
The calculations will of course come out similarly for a unit
current flow supplied across any other edge of $K_n$.  

The first distribution we are going to examine is of the degree
in $\Tree$ of a vertex, say $v_1$.  Since we are interested in
which of the edges incident to $v_1$ are in $\Tree$, we need
to calculate $H(\overline{v_1 v_i} , \overline{v_1 v_j})$ for
every $i,j \neq 1$.  Orienting all of these edges away from
$v_1$ and using the voltages we just worked out gives
$$H(\overline{v_1 v_i} , \overline{v_1 v_j}) = \left \{ \begin{array}{ccl} 
    ~~2/n~~ & : & ~~i = j \\ ~~1/n~~ & & \mbox{  otherwise} 
    \end{array} \right. . $$
Denoting the edge from $v_1$ to $v_i$ by $e_i$, we have the $n-1$ by
$n-1$ matrices
\vspace{2ex}
$$M(e_2 , \ldots , e_n) = \left ( \begin{array}{cccc} {2 \over n} & {1 
    \over n} & \cdots & {1 \over n}  \\ {1 \over n} & {2 \over n} & \cdots 
    & {1 \over n} \\ &&& \\ &&\vdots & \\
    &&& \\ {1 \over n} & {1 \over n} & \cdots & {2 \over n} \end{array} 
    \right ) \hspace{.4in}
 M^{(n-1)} (e_2 , \ldots , e_n) = \left ( \begin{array}{cccc} 
    {n-2 \over n} & {-1 \over n} & \cdots & {-1 \over n}  \\ {-1 \over n} 
    & {n-2 \over n} & \cdots & {-1 \over n} \\ &&& \\ &&\vdots & \\ &&& 
    \\ {-1 \over n} & {-1 \over n} & \cdots & {n-2 \over n} \end{array} 
    \right ). $$
\vspace{2ex}
There must be at least one edge in $\Tree$ incident to $v_1$ so
Theorem~\ref{incl-excl} says $\det M^{(n-1)} = \P (e_2 , \ldots , e_n
\notin \Tree) = 0$.  This is easy to verify: the rows sums to zero.
We can use $M^{(n-1)}$ to calculate the probability that $e_2$
is the only edge in $\Tree$ incident to $v_1$ by noting that this
happens if and only if $e_3 , \ldots , e_n \notin \Tree$.  This
is the determinant of $M^{(n-2)} (e_3 , \ldots , e_n)$ which is 
a matrix smaller by one thatn $M^{(n-1)}(e_2 , \ldots , e_n)$ but
which still has $(n-2)/n$'s down the diagonal and $-1/n$'s elsewhere.
This is a special case of a {\em circulant} matrix, which is a type
of matrix whose determinant is fairly easy to calculate.  

A $k$ by $k$ circulant matrix is an $M$ for which $M(i,j)$ 
is some number $a(i-j)$ depending only on $i-j$ mod $k$.  Thus
$M$ has $a_0$ all down the diagonal for some $a_0$, $a_1$ on the
next diagonal, and so forth.  The eigenvalues of a circulant
matrix $\lambda_0 , \ldots , \lambda_{k-1}$ are given by
$\lambda_j = \sum_{t=0}^{k-1} a_t \zeta^{jt}$ where $\zeta
= e^{2 \pi i/n}$ is the $n^{th}$ root of unity.  It is easy
to verify that these are the eigenvalues, by checking that
the vector $\ww$ for which $w_t = \zeta^{tj}$ is an eigenvector
for $M$ (no matter what the $a_i$ are) and has eigenvalue $\lambda_j$.
The determinant is then the product of the eigenvalues.
Details of this may be found in \cite{St}.  
 
In the case of $M^{(n-2)}$, $a_0 = (n-2)/n$ and $a_j = -1/n$
for $j \neq 0$.  Then $\lambda_0 = \sum_j a_j = 1/n$.  To calculate
the other eigenvalues note that for any $j \neq 0$ mod $n-2$, 
$\sum_{t=1}^{n-3} \zeta^{jt} = 0$.  Then $\lambda_j = (n-2)/n
\sum_{t=1}^{n-3} (-1/n) \zeta^{jt} = (n-1)/n - (1/n) \sum_{t=0}^{n-3}
\zeta^{tj} = (n-1)/n$.  This gives
$$\det M^{(n-2)} = \prod_{j=0}^{n-3} \lambda_j = {1 \over n} \,
  \left ( {n-1 \over n} \right )^{n-3} = {1+ o(1) \over ne}$$
as $n \rightarrow \infty$.  \footnote{Here, $o(1)$ signifies a quantity
going to zero as $n \rightarrow \infty$.  This is a convenient and 
standard notation that allows manipulation such as $(2+o(1))(3+o(1))=6+o(1)$.}  
Part of the Poisson limit has emerged: the probability
that $v_1$ has degree one in $\Tree$ is (by symmetry) $n-1$ times 
the probability that the particular edge $e_2$ is the only edge in 
$\Tree$ incident to $v_1$; this is $(n-1) (1+o(1))/en$ so it
converges to $e^{-1}$ as $n \rightarrow \infty$.  This is
$\P (X=1)$ where $X$ is one plus a Poisson(1) , i.e. a Poisson of mean
one. 

Each further part of the Poisson limit requires a more careful
evaluation of the limit.  To illustrate, we carry out the
second step.  Use one more degree of precision in the
Taylor series for $\ln (x)$ and $\exp (x)$ to get
\begin{eqnarray*}
&& n^{-1} \left ( {n-1 \over n} \right )^{n-3} \\[2ex]
& = & n^{-1} \exp [ (n-3) (-n^{-1} - n^{-2} (1/2 + o(1)))] \\[2ex]
& = & n^{-1} \exp [-1 + (5/2 + o(1)) n^{-1}] \\[2ex]
& = & n^{-1} e^{-1} [1 + (5/2 + o(1)) n^{-1}] .
\end{eqnarray*} 
The reason we need this precision is that we are going to calculate
the probability of $v_1$ having degree $2$ by summing the
$\P (e,f$ are the only edges incident to $v_1$ in $\Tree)$ over
all pairs of edges $e,f$ coming out of $v_1$.  By symmetry
this is just $(n-1)(n-2)/2$ times the probability that the particular
edges $e_2$ and $e_3$ are the only edges in $\Tree$ incident to $v_1$.
This probability is the determinant of a matrix which is not
a circulant, and to avoid calculating a difficult determinant 
it is better to write this probability as the following difference:
the probability that no edges other that $e_2$ and $e_3$ are incident
to $v_1$ minus the probability that $e_2$ is the only edge
incident to $v_1$ minus the probability that $e_2$ is the only edge
incident to $v_3$.  Since the final probability is this difference
multiplied by $(n-1)(n-2)/2$, the difference should be of order
$n^{-2}$, which explains why this degree of precision is required
for the latter two probabilities.  

The probability of $\Tree$ containing no edges incident to $v_1$ 
other than $e_2$ and $e_3$ is the determinant of $M^{(n-3)}(e_4 , 
\ldots ,e_n)$, which is an $n-3$ by $n-3$ circulant again having
$(n-2)/n$ on the diagonal and $-1/n$ elsewhere.  Then $\lambda_0
= \sum_{j=0}^{n-4} a_j = 2/n$ and $\lambda_j = (n-1)/n$ for 
$j \neq 0$ mod $n-3$, yielding 
$$\det M^{(n-3)} = 2 n^{-1} \left ( {n-1 \over n} \right )^{n-4}
= 2 n^{-1} e^{-1} [1 + (7/2 + o(1)) n^{-1}] $$
in the same manner as before.  Subtracting off the probabilities
of $e_2$ or $e_3$ being the only edge in $\Tree$ incident to $v_1$
gives
\begin{eqnarray*}
&&\P (e_2 , e_3 \in \Tree , e_4 , \ldots , e_n \notin \Tree) \\[2ex]
& = & 2n^{-1} e^{-1} [1+(7/2 +o(1))n^{-1}] - 2n^{-1} e^{-1} [1+(5/2
    +o(1))n^{-1}] = (2+o(1)) n^{-2} e^{-1} .
\end{eqnarray*}
Multiplying by $(n-1)(n-2)/2$ gives
$$\P (v_1 \mbox{ has degree 2 in } \Tree) \rightarrow e^{-1} $$
as $n \rightarrow \infty$, which is $\P (X=2)$ where $X$ is
one plus a Poisson(1).

\subsection{Another point of view} \label{5.2}

The calculations of the last section may be continued {\em ad infinitum}, 
but each 
step requires a more careful estimate so it pays to look for a way
to do all the steps at once.  The right alternative method will
be more readily apparent if we generalize to graphs
other than $K_n$ which do not admit such a precise calculation
(if a tool that is difficult to use breaks, you
may discover a better one).  

The important feature about $K_n$ was that the voltages were
easy to calculate.  There is a large class of graphs for which
the voltages are just as easy to calculate approximately.  
The term ``approximately'' can be made more rigorous by
considering sequences of graphs $G_n$ and stating approximations
in terms of limits as $n \rightarrow \infty$.  Since I've always
wanted to name a technical term after my dog, call a 
sequence of graphs $G_n$ {\em Gino-regular} if there is a sequence
$D_n$ such that
\begin{quote}
$(i)$  The maximum and minimum degree of a vertex in $G_n$ are $(1+o(1))D_n$
as $n \rightarrow \infty$; and \\

$(ii)$  The maximum and minimum over vertices $x \neq y,z$ of $G_n$
of the probability that $SRW_x^{G_n}$ hits $y$
before $z$ are $1/2 + o(1)$ as $n \rightarrow \infty$.
\end{quote}
Condition $(ii)$ implies that $D_n \rightarrow \infty$, so the
graphs $G_n$ are growing locally.
It is not hard to see that the voltage $V(z)$ in a unit current
flow across any edge $e=\xy$ of a graph $G_n$ in a Gino-regular
sequence is $(1+o(1)) D_n^{-1}  (\dd_x - \dd_y)(z)$ uniformly over
all choices of $x,y,z \in G_n$ as $n \rightarrow \infty$.  
The complete graphs $K_n$ are Gino-regular.  So are the $n$-cubes,
$B_n$, whose vertex sets are all the $n$-long sequences of zeros
and ones and whose edges connect sequences differing in only one
place.

\begin{picture}(200,100)
\put(20,20){\circle{5}}
\put(20,20){\line(1,0){40}}
\put(20,20){\line(0,1){40}}
\put(20,60){\circle{5}}
\put(20,60){\line(1,0){40}}
\put(60,20){\circle{5}}
\put(60,20){\line(0,1){40}}
\put(60,60){\circle{5}}
\put(33,5){$B_2$}
\put(110,20){\circle{5}}
\put(110,20){\line(1,0){40}}
\put(110,20){\line(0,1){40}}
\put(110,60){\circle{5}}
\put(110,60){\line(1,0){40}}
\put(150,20){\circle{5}}
\put(150,20){\line(0,1){40}}
\put(150,60){\circle{5}}
\put(123,5){$B_3$}
\put(130,40){\circle{5}}
\put(130,40){\line(1,0){40}}
\put(130,40){\line(0,1){40}}
\put(130,80){\circle{5}}
\put(130,80){\line(1,0){40}}
\put(170,40){\circle{5}}
\put(170,40){\line(0,1){40}}
\put(170,80){\circle{5}}
\put(110,20){\line(1,1){20}}
\put(110,60){\line(1,1){20}}
\put(150,20){\line(1,1){20}}
\put(150,60){\line(1,1){20}}
\end{picture}
figure~\thepfigure
\label{pfig5.2}
\addtocounter{pfigure}{1}

To see why $\{ B_n \}$ is Gino-regular, consider the ``worst case''
when $x$ is a neighbor of $y$.  There is a small probability that 
$SRW_x (1)$ will equal $y$, small because this is $\mbox{degree}(x)^{-1}
= (1+o(1)) D_n^{-1}$ which is going to zero.  There are even smaller
probabilities of reaching $y$ in the next few steps; in general,
unless $SRW_x$ hits $y$ in one step, it tends to get ``lost'' and
by the time it comes near $y$ or $z$ again it is thoroughly random
and is equally likely to hit $y$ or $z$ first.  In fact Gino-regular
sequences may be thought of as graphs that are nearly degree-regular,
which $SRW$ gets lost quickly.  

The approximate voltages give approximate transfer-impedances
$H(e,f) = (2+o(1))/n$ if $e = f$, $(1+o(1))/n$ if $e$ and $f$
meet at a single vertex (choose orientations away from the vertex)
and $o(1)/n$ if $e$ and $f$ do not meet.  The determinant of a
matrix is continuous in its entries, so it may seem that we have 
everything necessary to calculate limiting probabilities as
limits of determinants of transfer-impedance matrices.  If $v$ is a 
vertex in $G_k$ and $e_1 , \ldots , e_n$ are the edges incident
to $v$ in $G_k$ (so $n \approx D_k$), then the probability
of $e_2$ being the only edge in $\Tree$ incident to $v$ is
the determinant of 
$$M^{(n-1)} (e_2 , \ldots , e_n) = \left ( \begin{array}{cccc} (n-2+o(1))/n & 
    (-1+o(1))/n & \cdots & (-1+o(1))/n  \\ (-1+o(1))/n & (n-2+o(1))/n & 
    \cdots & (-1+o(1))/n \\ &&& \\ &&\vdots & \\ &&& \\ (-1+o(1))/n & 
    (-1+o(1))/n & \cdots & (n-2+o(1))/n \end{array} \right ). $$
Unfortunately, the matrix is changing size as $n \rightarrow \infty$,
so convergence of each entry to a known limit does not give us
the limit of the determinant.  

If the matrix were staying the same size, the problem would disappear.
This means we can successfully take the limit of probabilities of
events as long as they involve a bounded number of edges.  Thus
for any fixed edge $e_1$, $\P (e_1 \in \Tree) = \det M(e_1)
= (1+o(1)) (2/n)$.  For any fixed pair of edges $e_1$ and $e_2$ incident
to the same vertex,
$$\P (e_1 , e_2 \in \Tree)\; =\; \det M(e_1,e_2) = \left | \begin{array}{cc}
   (2+o(1))/n & (1+o(1)) / n \\ (1+o(1))/n & (2+o(1)) / n
   \end{array} \right | \;=\; (3+o(1)) n^{-2} .$$
In general if $e_1 , \ldots , e_r$ are all incident to $v$ then
the transfer-impedance matrix is $n^{-1}$ times an $r$ by $r$ matrix 
converging to the matrix with $2$ down the diagonal and $1$ elsewhere.  The
eigenvalues of this circulant are $\lambda_0 = r+1$ and $\lambda_j 
= 1$ for $j \neq 0$, yielding
$$\P (e_1 , \ldots , e_r \in \Tree) = (r+1+o(1)) n^{-r} .$$

What can we do with these probabilities?  Inclusion-exclusion fails
for the same reason as the large determinants fail -- the $o(1)$
errors pile up.  On the other hand, these probabilities determine
certain expectations.  Write $e_1 , \ldots , e_n$
again for the edges adjacent to $v$ and $I_i$ for
the indicator function which is one when $e_i \in \Tree$ and zero
otherwise; then   
$$\sum_i \P (e_i \in \Tree) = \sum_i \E I_i = \E \sum_i I_i = \E \deg (v) .$$
This tells us that $\E \deg (v) = n (2+o(1))n^{-1} = 2+o(1)$.
If try this with ordered pairs of edges, we get
$$\sum_{i \neq j} \P (e_i , e_j \in \Tree) = \sum_{i \neq j} \E I_i I_j = 
\E \sum_{i \neq j} I_i I_j .$$
This last quantity is the sum of all distinct ordered pairs of edges
incident to $v$ of the quantity: $1$ if they are both in the tree and 
0 otherwise.  If $\deg (v) = r$ then a one occurs in this sum $r(r-1)$ 
times, so the sum is $\deg (v) (\deg (v)-1)$.  The determinant calculation 
gave $\P (e_i,e_j \in \Tree) = (3+o(1))n^{-2}$ for each $i,j$, so
$$\E [\deg (v) (\deg(v)-1)] = n(n-1) (3+o(1))n^{-2} = 3+o(1) .$$
In general, using ordered $r$-tuples of distinct edges gives
\begin{eqnarray*}
&& \E [\deg (v) (\deg(v) - 1) \cdots (\deg(v)-r+1)] \\
& = & n(n-1) \cdots (n-r+1) (r+1+o(1))n^{-r} \\
& = & r+1+o(1) .
\end{eqnarray*}
Use the notation $(A)_r$ to denote $A(A-1)\cdots (A-r+1)$ which is
called the $r^{th}$ lower factorial of $A$.  If $Y_n$ is the random
variable $\deg (v)$ then we have succinctly,
\begin{equation} \label{moment 1}
\E (Y_n)_r = r+1+o(1) .
\end{equation}
$\E (Y_n)_r$ is called the $r^{th}$ factorial moment of $Y_n$.

If you remember why we are doing these calculations, you have probably
guessed that $\E (X)_r = r+1$ when $X$ is one plus a Poisson(1). 
This is indeed true and can be seen easily enough from the
logarithmic moment generating function $\E t^X$ via the identity
$$\E (X)_r = \left. \left ({d \over dt} \right )^r \right |_{t=1} \E t^X , $$
using $\E t^X = \E e^{X \ln (t)} = \phi (\ln (t)) = t e^{t-1}$ ; consult 
\cite[page 301]{Ro} for details.  All that we need now for a 
Poisson limit result is a theorem saying that
if the factorial moments of $Y_n$ are each converging to the 
factorial moments of $X$, then $Y_n$ is actually converging in
distribution to $X$.  This is worth spending a short subsection
on because it is algebraically very neat.

\subsection{The method of moments} \label{5.3}

A standard piece of real analysis shows that if all the factorial
moments of a sequence of random variables converging to a limit
are finite, then for each $r$, the limit of the $r^{th}$ factorial 
moments is the $r^{th}$ factorial moment of the limit.  
(This is essentially the Lebesgue-dominated convergence theorem.)
Another standard result is that if the moments of a sequence
of random variables converge, then the sequence, or at least
some subsequence is converging in distribution to some other
random variable whose moments are the limits of the moments in the
sequence.  Piecing together these straight-forward
facts leaves a serious gap in our prospective proof: What 
if there is some random variable $Z$ distributed differently from
$X$ with the same factorial moments?  If this could happen, then
there would be no reason to think that $Y_n$ converged in distribution
to $X$ rather than $Z$.  This scenario can actually happen -- there
really are differently distributed random variables with the same moments! 
(See the discussion of the lognormal distribution in \cite{Fe}.)  
Luckily this only happens when $X$ is badly behaved, and a
Poisson plus one is not badly behaved.  Here then is a proof
of the fact that the distribution of $X$ is the only one
with $r^{th}$ factorial moment $r+1$ for all $r$.  I will leave
it to you to piece together, look up in \cite{Fe} or take on faith how this
fact plus the results from real analysis imply $Y \dconv X$.

\begin{th} \label{moment method}
Let $X$ be a random variable with $\E (X)_r \leq e^{kr}$ for some $k$.
Then no random variable distributed differently from $X$ has the
same factorial moments.
\end{th}

\noindent{Proof:}  The factorial moments $\E (X)_r$ determine
the regular moments $\mu_r = \E X^r$ and {\em vice versa} by the
linear relations $(X)_1 = X^1 ; (X)_2 = X^2 - X^1$, etc.
 From these linear relations it also follows that factorial
moments are bounded by some $e^{kr}$ if and only if regular
moments are bounded by some $e^{kr}$, thus it suffices to
prove the theorem for regular moments.  Not only do the
moments determine the distribution, it is even possible to
calculate $\P (X=j)$ directly from the moments of $X$ in
the following manner.  

The {\em characteristic function} of $X$ is the function $\phi (t) = 
\E e^{itX}$ where $i = \sqrt{-1}$.  This is determined by the
moments since $\E e^{itX} = \E (1 + (itX) + (itX)^2/2! + \cdots )
= 1 + it\mu_1 + (it)^2 \mu_2 / 2! + \cdots$.  We use the exponential
bound on the growth of $\mu_r$ to deduce that this is absolutely
convergent for all $t$ (though a somewhat weaker condition would do).
The growth condition also shows that $\E e^{itX}$ is bounded and
absolutely convergent for
$y \in [0,2\pi]$.  Now $\P (X=j)$ can be determined by Fourier
inversion:
\begin{eqnarray*}
&& {1 \over 2 \pi} \int_0^{2\pi} \E e^{itX} e^{-ijt} dt  \\[2ex]
& = & {1 \over 2 \pi} \int_0^{2\pi} [\sum_{r \geq 0} e^{itr} \P (X=r)]
     e^{-ijt} dt  \\[2ex]
& = & {1 \over 2 \pi} \sum_{r \geq 0} \P (X=r) \int_0^{2\pi} e^{itr}
     e^{-ijt} dt  \\[2ex]
&& \mbox{ (switching the sum and integral is OK for bounded, absolutely
convergent integrals)} \\[2ex]
& = & {1 \over 2 \pi} \sum_{r \geq 0} \P (X=r) \dd_0 (r-j) \\[2ex]
& = & \P (X=j) .
\end{eqnarray*}
$\Cox$

\subsection{A branching process} \label{5.4}

In the last half of section~\ref{1.4} I promised to explain how 
convergence in distribution of $\deg (v)$ was a special case of
convergence of $\Tree$ near $v$ to a distribution called $\pois$.
(You might want to go back and reread that section before continuing.)
The infinite tree $\pois$ is interesting in its own right and 
I'll start making good on the promise by describing $\pois$.  

This begins with a short description of {\em Galton-Watson
branching processes}.  You can think of a Galton-Watson process
as a family tree for some fictional amoebas.  These fictional amoebas 
reproduce by splitting into any number of smaller amoebas (unlike
real amoebas that can only split into two parts at a time).  At
time $t=0$ there is just a single amoeba, and at each time 
$t=1,2,3,\ldots$, each living amoeba $\A$ splits into a random
number $N = N_t(\A)$ of amoebas, where the random numbers are independent
and all have the same distribution $\P (N_t (\A) = j) = p_j$.
Allow the possibility that $N = 0$ (the amoeba died) or that 
$N=1$ (the amoeba didn't do anything).  Let $\mu = \sum_j j p_j$
be the mean number of amoebas produced in a split.  A standard
result from the theory of branching processes \cite{AN} is that
if $\mu > 1$ then there is a positive probability that the
family tree will survive forever, the population exploding
exponentially as in the usual Malthusian forecasts for human
population in the twenty-first century.
Conversely when $\mu < 1$, the amoeba population dies out with 
probability 1 and in fact the chance of it surviving 
$n$ generations decreases exponentially with $n$.  When $\mu = 1$
the branching process is said to be {\em critical}.  It must
still die out, but the probability of it surviving $n$ generations
decays more slowly, like a constant times $1/n$.  The theory
of branching processes is quite large and you can find more 
details in \cite{AN} or \cite{Ha}.

Specialize now to the case where the random number of offspring has
a Poisson(1) distribution, i.e. $p_j = e^{-1} / j!$.
Here's the motivation for considering this case.  Imagine a graph 
$G$ in which each vertex has $N$ neighbors and $N$ is so large
it is virtually infinite.  Choose a subgraph $U$ by letting each edge
be included independently with probability $N^{-1}$.  Fix a vertex
$v \in G$ and look at the vertices connected to $v$ in $U$.  
The number of neighbors of $v$ in $U$ has a Poisson(1) 
distribution by the standard characterization of a Poisson as
the limit of number of occurrences of rare events.  For each neighbor
$y$ of $v$ in $U$, there are $N-1$ edges out of $y$ other than the 
one to $v$, and the number of those in $U$ will again be Poisson(1)
(since $N \approx \infty$, subtracting one does not matter)
and continuing this way shows that the connected component of
$v$ in $U$ is distributed as a Galton-Watson process with 
Poisson(1) offspring.  

Of course $U$ is not distributed like a uniform spanning tree $\Tree$.
For one thing, $U$ may with probability $e^{-1}$ fail to have any 
edges out of $v$.  Even if this doesn't happen, the chance of
$U$ having more than $n$ vertices goes to zero as $n \rightarrow
\infty$ (a critical Galton-Watson process dies out) whereas
$\Tree$, being a spanning tree of an almost infinite graph, goes on 
as far as the eye can see.  The next hope is that $\Tree$ looks like
$U$ conditioned not to die out.  This should in fact seem plausible:
you can check that $U$ has no cycles near $v$ since virtually all of
the $N$ edges out of each neighbor of $v$ lead further away from $v$;
then a uniform spanning tree should be a random cycle-free graph $U$ 
that treats each edge as equally likely, conditioned on being connected.

The conditioning must be done carefully, since the probability
of $U$ living forever is zero, but it turns out fine if you 
condition on $U$ living for at least $n$ generations and take the
limit as $n \rightarrow \infty$.  The random infinite tree $\pois$ that
results is called the {\em incipient infinite cluster} at $v$,
so named by percolation theorists (people who study connectivity
properties of random graphs).  It turns out there is an alternate
description for the incipient infinite cluster.  Let $v = v_0 , v_1 , 
v_2 , \ldots $ be a single line of vertices with edges $\overline{v
v_1} , \overline{v_1 v_2} , \ldots$.  For each of the vertices
$v_i$ independently, make a separate independent
copy $U_i$ of the critical Poisson(1) branching process $U$
with $v_i$ as the root and paste it onto the line already there.
Then this collage has the same distribution as $\pois$.  This
fact is the ``whole tree'' version of the fact that a Poisson(1)
conditioned to be nonzero is distributed as one plus a 
Poisson(1) (you can recover this fact from the fact about $\pois$ by
looking just at the neighbors of $v$).

\subsection{Tree moments}

To prove that a uniform spanning tree $\Tree_n$ of $G_n$ converges
in distribution to $\pois$ when $G_n$ is Gino-regular, we generalize
factorial moments to trees.  Let $t$ be a finite tree rooted 
at some vertex $x$ and let $W$ be a tree rooted at $v$.  $W$ is allowed
to be infinite but it must be locally finite -- only finitely
many edges incident to any vertex.  Say that a map $f$ from the vertices
of $t$ to the vertices of $W$ is a {\em tree-map} if $f$ is one to one,
maps $x$ to $v$ and neighbors to neighbors.  Let $N(W;t)$ count the
number of tree-maps from $t$ into $W$.  For example in the following
picture, $N(W;t) = 4$, since C and D can map to H and I in either 
order with A mapping to E, and B can map to F or G.  
\begin{picture}(200,140)
\put(20,20){\circle*{3}} 
\put(14,20){\scriptsize C}
\put(40,50){\circle*{3}}
\put(35,50){\scriptsize A}
\put(60,20){\circle*{3}}
\put(54,20){\scriptsize D}
\put(60,80){\circle*{3}}
\put(55,80){$x$}
\put(80,50){\circle*{3}}
\put(75,50){\scriptsize B}
\put(20,20){\line(2,3){20}}
\put(40,50){\line(2,3){20}}
\put(80,50){\line(-2,3){20}}
\put(60,20){\line(-2,3){20}}
\put(65,115){$t$}
\put(120,20){\circle*{3}}
\put(115,20){\scriptsize H}
\put(160,20){\circle*{3}}
\put(155,20){\scriptsize I}
\put(190,20){\circle*{3}}
\put(185,20){\scriptsize J}
\put(140,50){\circle*{3}}
\put(135,50){\scriptsize E}
\put(160,50){\circle*{3}}
\put(155,50){\scriptsize F}
\put(180,50){\circle*{3}}
\put(185,50){\scriptsize G}
\put(160,80){\circle*{3}}
\put(155,80){$v$}
\put(120,20){\line(2,3){20}}
\put(140,50){\line(2,3){20}}
\put(160,20){\line(-2,3){20}}
\put(180,50){\line(-2,3){20}}
\put(160,50){\line(0,1){30}}
\put(190,20){\line(-1,3){10}}
\put(155,115){$W$}
\end{picture}
figure~\thepfigure
\label{pfig5.3}
\addtocounter{pfigure}{1}

Define the $t^{th}$ {\em tree-moment} of a random tree $Z$ rooted
at $v$ to be $\E N(Z;t)$.  
If $t$ is an $n$-star, meaning a tree consisting of $n$ edges all
emanating from $x$, then a tree-map from $t$ to $W$ is just a
choice of $n$ distinct neighbors of $v$ in order, so $N(W;t) =
(\deg(v))_n$.  Thus $\E N(Z;t) = \E (\deg (v))_n$, the $n^{th}$
factorial moment of $\deg (v)$.  
This is to show you that tree-moments generalize
factorial moments.  Now let's see what the tree-moments of 
$\pois$ are.  Let $t$ be any finite tree and let $|t|$ denote 
the number of vertices in $t$.  

\begin{lem} \label{NUt=1}
Let $U$ be a Galton-Wastson process rooted at $v$ with 
Poisson(1) offspring.  Then $\E N(U;t) = 1$ for all finite trees $t$.
\end{lem}

\noindent{Proof:}  Use induction on $t$, the lemma being clear when
$t$ is a single vertex.  The way the induction step works for trees is 
to show that if a fact is true for a collection of trees $t_1 , \ldots
, t_n$ then it is true for the tree $t_*$ consisting of a root $x$ with
$n$ neighbors $x_1 , \ldots , x_n$ having subtrees $t_1 , \ldots t_n$ 
respectively as in the following illustration.

\begin{picture}(200,100)
\put(10,20){\circle*{3}}
\put(40,20){\circle*{3}}
\put(25,50){\circle*{3}}
\put(10,20){\line(1,2){15}}
\put(40,20){\line(-1,2){15}}
\put(25,5){$t_1$}
\put(65,20){\circle*{3}}
\put(64,5){$t_2$}
\put(90,20){\circle*{3}}
\put(80,50){\circle*{3}}
\put(90,20){\line(-1,3){10}}
\put(88,5){$t_3$}
\put(120,20){\circle*{3}}
\put(160,20){\circle*{3}}
\put(190,20){\circle*{3}}
\put(140,50){\circle*{3}}
\put(160,50){\circle*{3}}
\put(180,50){\circle*{3}}
\put(160,80){\circle*{3}}
\put(120,20){\line(2,3){20}}
\put(140,50){\line(2,3){20}}
\put(160,20){\line(-2,3){20}}
\put(180,50){\line(-2,3){20}}
\put(160,50){\line(0,1){30}}
\put(190,20){\line(-1,3){10}}
\put(160,5){$t_*$}
\end{picture}
figure~\thepfigure
\label{pfig5.4}
\addtocounter{pfigure}{1}

So let $t_1 , \ldots , t_n$ and $t_*$ be as above.  Any tree-map
$f : t_* \rightarrow U$ must map the $n$ neighbors of $v$ into 
distinct neighbors of $U$ and the expected number of ways to do
this is $\E (\deg (v))_n$ which is one for all $n$ since $\deg (v)$
is a Poisson(1) \cite{Fe}.  Now for any such assignment
of $f$ on the neighbors of $v$, the number of ways of completing
the assignment to a tree-map is the product over $i= 1 , \ldots n$
of the number of ways of mapping each $t_i$ into the subtree of 
$U$ below $f(x_i)$.  After conditioning on what the first generation
of $U$ looks like, the subtrees below any neighbors of $v$ are
independent and themselves Galton-Watsons with Poisson(1) offspring.
(This is what it means to be Galton-Watson.)  By induction then,
the expected number of ways of completing the assignment of $f$
is the product of a bunch of ones and is therefore one.  Thus
$\E N(U;t) = \E (\deg (v))_n \prod_{i=1}^n 1 = 1$.   $\Cox$

Back to calculating $\E N(\pois ; t)$.  Recall that $\pois$ is
a line $v_0 , v_1 , \ldots $ with Poisson(1) branching processes
$U_i$ stapled on.  Each tree-map $f : t \rightarrow \pois$
hits some initial segment $v_0 , \ldots v_k$ of the original
line, so there is some vertex $y_f \in t$ such that $f(y_f) = v_k$
for some $k$ but $v_{k+1}$ is not in the image of $f$.  For each
$y \in t$, we count the expected number of tree-maps $f$
for which $y_f = y$.  There is a path $x = f^{-1} (v_0) , \ldots , f^{-1} 
(v_k) = y$ in $t$ going from the root $x$ to $y$.
The remaining vertices of $t$ can be separated into $k+1$ subtrees
below each of the $f^{-1} (v_i)$.  These subtrees must then
get mapped respectively into the $U_i$.  By the lemma, the expected
number of ways of mapping anything into a $U_i$ is one, so the
expected number of $f$ for which $y_f = y$ is $\prod_{i=1}^k 1 = 1$.
Summing over $y$ then gives
\begin{equation}
\E N(\pois ; t) = |t| \end{equation}

The last thing we are going to do to in proving the stronger Poisson
convergence theorem is to show 
\begin{lem} \label{NTt}
Let $G_n$ be a Gino-regular sequence of graphs, and let $\Tree_n$
be a uniform spanning tree of $G_n$ rooted at some $v_n$.  Then
for any finite rooted tree $t$, $\E N(\Tree_n ; t) \rightarrow
|t|$ as $n \rightarrow \infty$.
\end{lem}
It is not trivial from here to establish that $\Tree \wedge r$
converges in distribution to $\pois \wedge r$ for every $r$.  
The standard real analysis facts I quoted in section~\ref{5.3}
about moments need to be replace by some not-so-standard (but
not too hard) facts about tree-moments.  Suffice it to say that
the previous two lemmas do in the end prove (see \cite{BP} for details)
\begin{th} \label{strong converge}
Let $G_n$ be a Gino-regular sequence of graphs, and let $\Tree_n$
be a uniform spanning tree of $G_n$ rooted at some $v_n$.  Then
for any $r$, $\Tree_n \wedge r$ converges in distribution to
$\pois \wedge r$ as $n \rightarrow \infty$.
\end{th}

\noindent{Sketch of proof of Lemma~\protect{\ref{NTt}}:}  Fix
a finite $t$ rooted at $x$.  To calculate the expected number
of tree-maps from $t$ into $\Tree_n$ we will sum over every possible
image of a tree-map the probability that all of those edges
are actually present in $\Tree_n$.  By an image of a tree-map,
I mean two things: (1) a collection $\{ v_x : x \in t \}$ of vertices of
$G_n$ indexed by the vertices of $t$ for which $v_x \sim v_y$
in $G$ whenever $x \sim y$ in $t$; (2) a collection of edges
$e_\ee$ connecting $v_x$ and $v_y$ for every edge $\ee \in t$ connecting
some $x$ and $y$.   Fix such an image.

The transfer-impedance theorem tells us that the probability
of finding all the edges $v_e$ in $\Tree$ is the determinant
of $M(e_\ee : \ee \in t)$.  Now for edges $e , e' \in G$,
Gino-regularity gives that $H(e , e') = D_n^{-1} (o(1) + \kappa)$
uniformly over edges of $G_n$,
where $\kappa$ is $2,1$ or $0$ according to whether $e=e'$,
they share an endpoint, or they are disjoint.  The determinant
is then well approximated by the corresponding determinant without
the $o(1)$ terms, which can be worked out as exactly $|t| D_n^{1-|t|}$.

This must now be summed over all possible images, which amounts
to multiplying $|t| D_n^{1-|t|}$ by the number of possible images.
I claim the number of possible images is approximately $D_n^{|t|-1}$.
To see this, imagine starting at the root $x$, which must get mapped
to $v_n$, and choosing successively where to map each nest vertex of
$t$.  Since there are approximately $D_n$ edges coming out of each
vertex of $G_n$, there are always about $D_n$ choices for the 
image of the next vertex (the fact that you are not allowed to
choose any vertex already chosen is insignificant as $D_n$ gets
large).  There are $|t|-1$ choices, so the number of maps
is about $D_n^{|t|-1}$.  This proves the claim.  The claim
implies that the expected number of tree-maps from $t$ to $\Tree_n$
is $|t| D_n^{1-|t|} D_n^{|t|-1} = |t|$, proving the lemma.    $\Cox$

\section{Infinite lattices, dimers and entropy}

There is, believe it or not, another model that ends up being equivalent
to the uniform spanning tree model under a correspondence at least as
surprising as the correspondence between spanning trees and random walks.
This is the so-called {\em dimer} or {\em domino tiling} model, which was
studied by statistical physicists quite independently
of the uniform spanning tree model.  The present section is 
intended to show how one of the fundamental questions of 
this model, namely calculating its entropy, can be solved using
what we know about spanning trees.  Since it's getting late,
there will be pictures but no detailed proofs.

\subsection{Dimers}

A dimer is a substance that on the molecular level is made up
of two smaller groups of atoms (imagine two spheres of matter)
adhering to each other via a covalent bond; consequently it is shaped
like a dumbbell.  If a bunch of dimer molecules are packed together
in a cold room and a few of the less significant laws of physics
are ignored, the molecules should array themselves into some
sort of regular lattice, fitting together as snugly as dumbbells can.
To model this, let $r$ be some positive real number representing the 
length of one of the dumbbells.  Let $L$ be a lattice, i.e. a 
regular array of points in three-space, for which each point in $L$
has some neighbors at distance $r$.  For example $r$ could be
$1$ and $L$ could be the standard integer lattice
$\{ (x,y,z) : x,y,z \in \Z \}$, so $r$ is the minimum distance
between any two points of $L$ (see the picture below).  
Alternatively $r$ could be $\sqrt{2}$
or $\sqrt{3}$ for the same $L$.  Make a graph $G$ whose vertices are
the points of $L$, with an edge between any pair of points
at distance $r$ from each other.  Then the possible packings
of dimers in the lattice are just the ways of partitioning the
lattice into pairs of vertices, each pair (representing one molecule) 
being the two enpoints of some edge.  The following picture shows
part of a packing of the integer lattice with nearest-neighbor edges.  \\
\begin{picture}(100,100)(-50,-10)
\put(0,0){\circle* {7}}
\put(0,30){\circle* {7}}
\put(0,60){\circle* {7}}
\put(30,30){\circle* {7}}
\put(30,60){\circle* {7}}
\put(20,50){\circle* {7}}
\put(50,50){\circle* {7}}
\put(20,80){\circle* {7}}
\put(50,80){\circle* {7}}
\put(60,30){\circle* {7}}
\multiput(-3,0)(.5,0){13}{\line(0,1){30}}
\multiput(30,27)(0,.5){13}{\line(1,0){30}}
\multiput(20,47)(0,.5){13}{\line(1,0){30}}
\multiput(-3,60)(.5,0){13}{\line(1,1){20}}
\multiput(27,60)(.5,0){13}{\line(1,1){20}}
\put(0,-10){\line(0,1){80}}
\put(30,20){\line(0,1){50}}
\put(20,40){\line(0,1){50}}
\put(50,40){\line(0,1){50}}
\put(-10,30){\line(1,0){80}}
\put(-10,60){\line(1,0){50}}
\put(20,50){\line(1,0){50}}
\put(20,80){\line(1,0){50}}
\put(-7,23){\line(1,1){34}}
\put(23,23){\line(1,1){34}}
\put(-7,53){\line(1,1){34}}
\put(23,53){\line(1,1){34}}
\end{picture}
figure~\thepfigure
\label{pfig6.1}
\addtocounter{pfigure}{1}

Take a large finite box inside the lattice, containing $N$ vertices.
If $N$ is even and the box is not an awkward shape, there
will be not only one but many ways to pack it with
dimers.  There will be several edges incident to each vertex $v$,
representing a choice to be made as to which other vertex 
will be covered by the molecule with one atom covering $v$.
These choices obviously cannot be made independently, but it
should be plausible from this that the total number of configurations
is approximately $\gamma^N$ for some $\gamma > 1$ as $N$
goes to infinity.  This number can be written alternatively
as $e^{hN}$ where $h = \ln (\gamma )$ is called the entropy
of the packing problem.  The thermodynamics of the resulting 
substance depend on, among other things, the entropy $h$.  

The case that has been studied the most is where $L$ is the 
two-dimensional integer lattice with $r=1$.  The graph $G$ is then
the usual nearest-neighbor square lattice.  Physically this
corresponds to packing the dimers between two slides.
You can get the same packing problem by attempting to tile
the plane with dominos -- vertical and horizontal 1 by 2 rectangles --
which is why the model also goes by the name of domino tiling.

\subsection{Dominos and spanning trees}

We have not yet talked about spanning trees of an infinite
graph, but the definition remains the same: a connected subgraph
touching each vertex and contaning no cycles.  If the subgraph
need not be connected, it is a spanning forest.  Define an
{\em essential spanning forest} or ESF to be a spanning forest
that has no finite components.  Informally, an ESF is a subgraph
that you can't distinguish from a spanning tree by only looking
at a finite part of it (since it has no cycles or {\em islands}).

Let $G_2$ denote the nearest-neighbor graph on the two dimensional
integer lattice.  Since $G_2$ is a planar graph, it has a {\em
dual} graph $G_2^*$, which has a vertex in each cell of $G_2$ and 
an edge $e^*$ crossing each edge $e$ of $G_2$.  In the following picture,
filled circles and heavy lines denote $G_2$ and open circles 
and dotted lines denote $G_2^*$.  Note that $G_2$, together with $G_2^*$
and the points where edges cross dual edges, forms another graph
$\tilde{G_2}$ that is just $G_2$ scaled down by a factor of two.   \\
\begin{picture}(160,100)(-65,0)
\put(10,20){\line(1,0){80}}
\put(10,50){\line(1,0){80}}
\put(10,80){\line(1,0){80}}
\put(20,10){\line(0,1){80}}
\put(50,10){\line(0,1){80}}
\put(80,10){\line(0,1){80}}
\put(20,20){\circle*{3}}
\put(20,50){\circle*{3}}
\put(20,80){\circle*{3}}
\put(50,20){\circle*{3}}
\put(50,50){\circle*{3}}
\put(50,80){\circle*{3}}
\put(80,20){\circle*{3}}
\put(80,50){\circle*{3}}
\put(80,80){\circle*{3}}
\put(35,35){\dashbox{2}(30,30)}
\put(35,35){\circle{3}}
\put(35,65){\circle{3}}
\put(65,35){\circle{3}}
\put(65,65){\circle{3}}
\put(5,35){\dashbox{2}(30,0)}
\put(5,65){\dashbox{2}(30,0)}
\put(65,35){\dashbox{2}(30,0)}
\put(65,65){\dashbox{2}(30,0)}
\put(35,5){\dashbox{2}(0,30)}
\put(65,5){\dashbox{2}(0,30)}
\put(35,65){\dashbox{2}(0,30)}
\put(65,65){\dashbox{2}(0,30)}
\end{picture}
figure~\thepfigure
\label{pfig6.2}
\addtocounter{pfigure}{1} \\
Each subgraph $H$ of $G$ has a dual subgraph $H^*$ consisting
of all edges $e^*$ of $G^*$ dual to edges $e$ \underline{not}
in $H$.  If $H$ has a cycle, then the duals of all edges in the
cycle are absent from $H^*$ which separates $H^*$ into two
components: the interior and exterior of the cycle.  Similarly,
an island in $H$ corresponds to a cycle in $H^*$ as in the picture:
\begin{picture}(180,100)(-50,10)
\multiput(20,20)(30,0){4}{\circle*{3}}
\multiput(20,50)(30,0){4}{\circle*{3}}
\multiput(20,80)(30,0){4}{\circle*{3}}
\multiput(35,35)(30,0){3}{\circle {3}}
\multiput(35,65)(30,0){3}{\circle {3}}
\put(20,20){\line(1,0){90}}
\put(20,20){\line(0,1){60}}
\put(110,20){\line(0,1){60}}
\put(50,80){\line(1,0){60}}
\put(50,50){\line(1,0){30}}
\put(20,20){\line(0,1){60}}
\put(35,35){\dashbox{2}(60,30)}
\put(35,65){\dashbox{2}(0,30)}
\end{picture}
figure~\thepfigure
\label{pfig6.3}
\addtocounter{pfigure}{1} \\
 From this description, it is clear that $T$ is an essential spanning
forest of $G_2$ if and only if $T^*$ is an essential spanning forest 
of $G_2^* \cong G_2$.  

Let $T$ now an infinite tree.  We
define {\em directed} a little differently than in the finite case:
say $T$ is directed if the edges are oriented so that every
vertex has precisely one edge leading out of it.  Following the
arrows from any vertex gives an infinite path and it is not hard
to check that any two such paths from different vertices eventually
merge.  Thus directedness for infinite trees is like directedness
for finite trees, toward a vertex at infinity.  

Say an essential spanning forest of $G_2$ is directed
if a direction has been chosen for each of its components and each
of the components of its dual.  Here then is the connection between
dominos and essential spanning forests.  
\begin{quotation}
Let $T$ be a directed essential spanning forest of $G_2$, with
dual $T^*$.  Construct a domino tiling of $\tilde{G_2}$ as follows.
Each vertex $v \in V(G_2) \subseteq V(\tilde{G_2})$ is covered by a 
domino that also covers the vertex of $\tilde{G_2}$ in the middle
of the edge of $T$ that leads out of $v$.  Similarly,
each vertex $v^* \in V(G_2^*)$ is covered by a 
domino also covering the middle of the edge of $T^*$ leading out
of $v$.  It is easy to check that this gives a legitimate domino
tiling: every domino covers two neighboring vertices,
and each vertex is covered by precisely one domino.  \\

Conversely, for any domino tiling of $\tilde{G_2}$, directed essential
spanning forests $T$ and $T^*$ for $G_2$ and $G_2^*$ can be constructed
as follows.  For each $v \in V(G_2)$, the oriented edge leading
out of $v$ in $T$ is the one along which the domino covering $v$
lies (i.e. the one whose midpoint is the other vertex of $\tilde{G_2}$
covered by the domino covering $v$).  Construct $T^*$ analogously. 
To show that $T$ and $T^*$ are directed ESF's amounts to showing
there are no cycles, since clearly $T$ and $T^*$ will have one edge 
coming out of each vertex.  This is true because if you set up
dominos in such a way as to create a cycle, they will always enclose an
odd number of vertices (check it yourself!).   Then there is no way to 
extend this configuration to a legitimate domino tiling of $\tilde{G_2}$.
\end{quotation}

It is easy to see that the two operations above invert each other,
giving a one to one correspondence between domino tilings of
$\tilde{G_2}$ and directed essential spanning forests of $G_2$.  
To bring this back into the realm of finite graphs requires 
ironing out some technicalities which I am instead going to
ignore.  The basic idea is that domino tilings of the $2n$-torus
$T_{2n}$ correspond to spanning trees of $T_n$ almost as well
as domino tilings of $\tilde{G_2}$ correspond to spanning trees of $G_2$.
Going from directed essential spanning forests to spanning trees
is one of the details glossed over here, but explained somewhat in the
next subsection.  The entropy for domino tilings is then one
quarter the entropy for spanning trees, since $T_{2n}$ has
four times as many vertices as $T_n$.  Entropy for spanning trees
just means the number $h$ for which $T_n$ has approximately
$e^{hn^2}$ spanning trees.  To calculate this, we use the
matrix-tree theorem.

The number of spanning trees of $T_n$ according to this theorem
is the determinant of a minor of the matrix indexed by vertices of $T_n$
whose $v,w$-entry is $4$ if $v=w$, $-1$ if $v \sim w$ and $0$
otherwise.  If $T_n$ were replaced by $n$ edges in a circle,
then this would be a circulant matrix.  As is, it is a 
generalized circulant, with symmetry group $T_n = (\Z / n\Z)^2$ 
instead of $Z / n\Z$.  The eigenvalues can be gotten via
group representations of $T_n$, resulting in eigenvalues
$4 - 2 \cos (2 \pi k/n) - 2 \cos (2 \pi l/n)$ as $k$ and
$l$ range from $0$ to $n-1$.  The determinant we want
is the product of all of these except for the zero eigenvalue at
$k=l=0$.  The log of the determinant divided by $n^2$
is the average of these as $k$ and $l$ vary, and the entropy
is the limit of this as $n \rightarrow \infty$ which is given by
$$\int_0^1 \int_0^1 \ln (4 - 2 \cos (2 \pi x) - 2 \cos (2 \pi y))
   \,\, dx \, dy . $$

\subsection{Miscellany}

The limit theorems in Section~\ref{5} involved letting $G_n$
tend to infinity locally, in the sense that each vertex
in $G_n$ had higher degree as $n$ grew larger.  Instead, one
may consider a sequence such as $G_n = T_n$; clearly the $n$-torus
converges in some sense to $G_2$ as $n \rightarrow \infty$, so there
ought to be some limit theorem.  Let $\Tree_n$ be a uniform
spanning tree of $G_n$.  Since $G_n$ is not Gino-regular, the
limit may not be $\pois$ and in fact cannot be since the limit
has degree bounded by four.  It turns out that $\Tree_n$ converges
in distribution to a random tree $\Tree$ called the uniform
random spanning tree for the integer lattice.  This works also
for any sequence of graphs converging to the three or four
dimensional integer lattices \cite{Pe}.  Unfortunately the process breaks
down in dimensions five and higher.  There the uniform spanning
spanning trees on $G_n$ do converge to a limiting distribution
but instead of a spanning tree of the lattice, you get an
essential spanning forest that has infinitely many components.
If you can't see how the limit of spanning trees could be a spanning
forest, remember that an essential spanning forest is so similar
to a spanning tree that you can't tell them apart with any finite
amount of information.  

Another result from this study is that in dimensions $2,3$ and $4$,
the uniform random spanning tree $\Tree$ has only one path
to infinity.  What this really means is that any two infinite
paths must eventually join up.  Not only that, but $\Tree^*$
has the same property.  That means there is only one way to
direct $\Tree$, so that each choice of $\Tree$ uniquely 
determines a domino tiling of $\tilde{G_2}$.  In this way
it makes sense to speak of a uniform random domino tiling 
of the plane: just choose a uniform random spanning tree
and see what domino tiling it corresponds to.  

That takes care of one of the details glossed over in the
previous subsection.  It also just about wraps up what I wanted
to talk about in this article.  As a parting note, let me mention
an open problem.  Let $G$ be the infinite nearest neighbor
graph on the integer lattice in $d$ dimensions and let $\Tree$
be the uniform spanning tree on $G$ gotten by taking a distributional
limit of uniform spanning trees on $d$-dimensional $n$-tori
as $n \rightarrow \infty$ as explained above.  

\begin{conj} \label{one end}
Suppose $d \geq 5$.  Then with probability one, each component 
of the essential spanning forest has only one path to infinity,
in the sense that any two infinite paths must eventually merge.
\end{conj}

\renewcommand{\baselinestretch}{1.0}\large
\end{document}